\documentclass[a4paper,12pt]{amsart}

\usepackage{amssymb}
\usepackage{xy}
\usepackage{bm}
\xyoption{all}
\xyoption{poly}

\newtheorem{theo}{Theorem}[section]
\newtheorem{coro}[theo]{Corollary}
\newtheorem{prop}[theo]{Proposition}
\newtheorem{lemma}[theo]{Lemma}

\theoremstyle{definition}
\newtheorem{defi}[theo]{Definition}
\newtheorem{remark}[theo]{Remark}
\newtheorem{notation}[theo]{Notation}
\newtheorem{example}[theo]{Example}
\newtheorem{convention}[theo]{Convention}

\def\CC{{\mathcal{C}}}
\def\CS{{\mathcal{S}}}

\def\CN{{\mathcal{N}}}

\def\CA{{\mathcal{A}}}

\def\CG{{\mathcal{G}}}

\def\CS{{\mathcal{S}}}

\def\CX{{\mathcal{X}}}

\def\BC{{\mathbf{C}}}

\def\BC{{\mathbb{C}}}

\def\BZ{{\mathbb{Z}}}

\def\op{\operatorname{op}\nolimits}

\def\GL{\operatorname{GL}\nolimits}

\def\Hom{\operatorname{Hom}\nolimits}

\def\Ext{\operatorname{Ext}\nolimits}

\def\reg{\operatorname{reg}\nolimits}

\def\End{\operatorname{End}\nolimits}

\def\Flag{\operatorname{Flag}\nolimits}

\def\bar{\operatorname{bar}\nolimits}
\def\gar{\operatorname{gar}\nolimits}
\def\Gal{\operatorname{Gal}\nolimits}

\def\iim{\operatorname{Im}\nolimits}
\def\ie{{\em i.e.}}

\title{Garside categories, periodic loops and cyclic sets}
\author{David Bessis}
\thanks{First draft,
some details and proofs are missing.}
\address{DMA, \'Ecole normale sup\'erieure, 
45 rue d'Ulm, 75230 Paris cedex 05, France}
\email{david dot bessis at ens dot fr}

\begin{document}

\begin{abstract}
Garside groupoids, as recently introduced by Krammer, 
generalise Garside groups.
A \emph{weak Garside group} is a group that is equivalent as a category
to a Garside groupoid. 
We show that any periodic loop in a Garside groupoid $\CG$ may be viewed
as a Garside element for a certain Garside
structure on another Garside groupoid $\CG_m$, which is equivalent as
a category to $\CG$.
As a consequence, the centraliser of a periodic element
in a weak Garside group is a weak Garside group.
Our main tool is the notion of \emph{divided Garside categories},
an analog for Garside categories of B\"okstedt-Hsiang-Madsen's subdivisions
of Connes' cyclic category.
This tool is used in our separate proof of the $K(\pi,1)$ property
for complex reflection arrangements.
\end{abstract}

\maketitle

\tableofcontents

\section*{Introduction}

A classical theorem says that a periodic homeomorphism of
the disk is conjugate to a rotation.
It was simultaneously announced by
Ker\'ekj\'art\'o and Brouwer in 1919 and the
first undisputed proof was published by Eilenberg in 1934
(see \cite{brouwer,eilenberg,kerekjarto}).
Our results include some analogs of 
Ker\'ekj\'art\'o-Brouwer-Eilenberg's theorem, in the context
of braid groups of complex reflection groups and, more generally,
cyclic Garside groupoids.

Let $X$ be a $S^1$-space, \ie, a topological space
together with a continuous action of
$S^1:=\{z\in\BC | |z|=1\}$.
The fundamental group $\pi_1(X,x_0)$ admits a special element,
called \emph{full-twist at $x_0$} and denoted by $\tau_{x_0}$
(or simply $\tau$),
represented by the path
$[0,1]\to X, t\mapsto e^{2i\pi t}x_0$. It is clearly central.
The map $x_0\mapsto \tau_{x_0}$ lies in the ``centre'' of the
fundamental groupoid, in the sense that it is a natural automorphism
of the identity functor.

More generally, an element
of $\pi_1(X,x_0)$ is a \emph{rotation} of angle $\theta$ if
it is represented by $t\mapsto e^{i\theta t}x_0$. This of course
requires $x_0$ to be $e^{i\theta}$-invariant.
It is enough to restrict one's attention to \emph{rational rotations},
those whose angles are rational multiples of $\pi$, because for other values
of $\theta$ the basepoint $x_0$ must be $S^1$-fixed, which implies that
irrational rotations are trivial.
The full-twist is a rotation
of angle $2\pi$. Rotations may be composed and their angles add up.

Let $p\in \BZ, q\in\BZ_{\geq 1}$.
An element $\gamma\in \pi_1(X,x_0)$ is \emph{$\frac{p}{q}$-periodic}
(or simply \emph{periodic}) if
$\gamma^q=\tau^p$ (note that this may \emph{a priori} depend on the actual
$p,q$ and not just on the rational number they represent). 
A rotation of angle $2\pi\frac{p}{q}$ is $\frac{p}{q}$-periodic.

Given any $S^1$-space $X$, one may be interested in the following questions,
whose answers obviously depend on $X$:

\begin{quote}
Question 1. \emph{
Is any periodic element conjugate to a rotation ?
}
\end{quote}

\begin{quote}
Question 2. \emph{
Are two rotations with identical angles conjugate ?
}
\end{quote}

In Question 2, ``conjugate'' should be understood
in the groupoid sense: $\rho\in\pi_1(X,x)$ and $\rho'\in\pi_1(X,x')$ are
conjugate if there exists $\gamma\in\Hom_{\pi_1(X)}(x,x')$ such that
$\rho\gamma=\gamma\rho'$.

There is an easy sufficient condition for Question 2 to have a positive
answer: assume that $p\wedge q=1$;
 let $\mu_q\subseteq S^1$ be the subgroup of $q$-th roots of unity
and let $X^{\mu_q}\subseteq X$ be the corresponding set of fixed points;
if $\pi_1(X^{\mu_q})$ is a connected category (or, equivalently,
is $X^{\mu_q}$ is path-connected), then any two rotations $\rho,\rho'$
of identical angles $2\pi\frac{p}{q}$ are conjugate.
Indeed, let $x$ be the basepoint of $\rho$, $x'$ be the basepoint of $\rho'$.
Because we have assume $p\wedge q=1$, $x$ and $x'$ are in $X^{\mu_q}$.
For any $\gamma$ in the direct image of $\Hom_{\pi_1(X^{\mu_q})}(x,x')$
under the natural functor $\pi_1(X^{\mu_q}) \to \pi_1(X)$,
there is an obvious homotopy showing that $\rho\gamma=\gamma\rho'$.

Actually, the restriction to $X^{\mu_q}$ of the $S^1$-structure on $X$
admits a ``$q$-th root''
\begin{eqnarray*}
[0,1] \times X^{\mu_q} & \longrightarrow & X^{\mu_q} \\
(t,x) & \longmapsto & e^{\frac{2i\pi t}{q}} x
\end{eqnarray*}
whose full-twists are preimages via $\pi_1(X^{\mu_q}) \to \pi_1(X)$
of rotations of angle $\frac{2\pi}{q}$.
As a consequence, rotations of angle $2\pi\frac{p}{q}$ lie
in the centre of $\iim(\pi_1(X^{\mu_q}) \to \pi_1(X))$.
This leads to the following questions:

\begin{quote}
Question 3 (group version). \emph{
Let $\rho$ be a rotation of angle $2\pi\frac{p}{q}$, with $p\wedge q=1$.
Let $x_0$ be the basepoint of $\rho$.
Does the inclusion $X^{\mu_q}\hookrightarrow X$ induce an isomorphism
$$\pi_1(X^{\mu_q},x_0)\simeq C_{\pi_1(X,x_0)}(\rho) \text{ ?}$$
}
\end{quote}

\begin{quote}
Question 3' (groupoid version). \emph{
Let $q$ be a positive integer, let $p\in \BZ$ such that $p\wedge q =1$.
Is the natural functor
$\pi_1(X^{\mu_q})\to \pi_1(X)$ faithful? Does it identify
$\pi_1(X^{\mu_q})$ with the largest subcategory of $\pi_1(X)$ with object
set $X^{\mu_q}$ and whose ``centre'' contains the natural family
of rotations of angle $2\pi\frac{p}{q}$?
}
\end{quote}

\begin{example}
\label{example01}
Let $X_n$ be the space of unordered configurations of $n$ distinct
points in $\BC$, whose fundamental group is the classical braid group $B_n$
on $n$ strings.
It has a natural structure of $S^1$-space (via the multiplicative action of
$S^1$ on $\BC$). 
The theorem of Ker\'ekj\'art\'o-Brouwer-Eilenberg implies
that Question 1 has a positive answer (to see this, one uses the
standard interpretation of the mapping class group of the punctured
disk as the quotient of $B_n$ by the full-twist). Because $X_n^{\mu_q}$
is either empty or connected, Question 2 has a positive answer.
Questions 3 and 3' also have positive answers, as it was shown in my
earlier work with Digne and Michel, \cite{bdm}.
\end{example}

\begin{example}
\label{myexample}
Let $V$ be a finite dimensional complex vector space. Let
$W\subseteq \GL(V)$ be a finite complex reflection group.
Let $V^{\reg}$ be the complement of the reflecting hyperplanes.
Scalar multiplication on $V$ commutes with $W$-action and preserves $V^{\reg}$.
It induces a structure of $S^1$-space on the quotient
$W\backslash V^{\reg}$, the
\emph{regular orbit space of $W$}. When $W$ is the symmetric group
in its permutation representation, the regular orbit space is $X_n$.
The fundamental group of $W\backslash V^{\reg}$ is the
\emph{(generalised) braid group} of $W$.
One application of the tools presented
here is the proof in \cite{kapiun} that, when $W$ is \emph{well-generated},
all above questions have positive answers when applied to
$W\backslash V^{\reg}$.
This theorem should be understood as a braid analog 
of Springer's theory of \emph{regular elements in complex reflection
groups}, \cite{springer}. 
It was motivated by a series of questions and conjectures by Brou\'e and Michel,
\cite{brmi,boston}, that naturally arose as they studied Deligne-Lusztig
varieties with the insight that their cohomology should provide tilting
modules verifying Brou\'e's abelian defect conjecture
for finite groups of Lie type.
In the particular case when $W$ is the symmetric group, we recover
Example \ref{example01}.
\end{example}

Our main result may be phrased in several ways, one of which is the
following:

\begin{theo}
\label{mykerekjarto}
If $X$ is the geometric realisation of the Garside nerve of
a cyclic Garside groupoid, then Question 1 has a positive answer.
\end{theo}

Moreover, in that setting, there are practical ways to tackle
Questions 2, 3 and 3'.

Before explaining what \emph{Garside nerves} and \emph{cyclic Garside 
groupoids} are, let us begin by discussing a particular case of Question 3
that admits an elementary solution relying on classical tools.

If $W$ is a complexified real reflection group,
it was proved by Brieskorn, \cite{brieskorn},
that the associated braid group
$B(W)$ is isomorphic to the  \emph{Artin group}
$$A(W,S):= \left< S \left| \underbrace{sts...}_{m_{s,t} \text{ terms}}
= \underbrace{tst...}_{m_{s,t} \text{ terms}} \right.\right>,$$
where $S$ is the set of reflections with respect to the walls of a chosen chamber of the real arrangement
and $(m_{s,t})$ is the Coxeter matrix.
Let $\Delta$ be the image by Tits' section $W\to A(W,S)$ of the longest
element $w_0$. One may choose the isomorphism $B(W)\simeq A(W,S)$ such
that $\Delta$ is a rotation of angle $\pi$, corresponding to the case
$p=1,q=2$ in Question 3. The word problem for $A(W,S)$ was solved
independently by Deligne and Brieskorn-Saito, \cite{deligne,brieskornsaito}.
In modern terms, their solution
relies on the fact that $A(W,S)$ is a \emph{Garside group}
with \emph{Garside element} $\Delta$ (see Section \ref{sectiongarside}
for more details about Garside groups).
Because the conjugacy action of $\Delta$ on $A(W,S)$ can be understood
through the Garside normal form, the centraliser is easy to compute.
Indeed, the conjugacy action of $w_0$ on $W$ is a 
\emph{diagram automorphism} (\ie, it is induced by a permutation of 
$S$) and the centraliser $W':=C_{W}(w_0)$ is a Coxeter group
with Coxeter generating set $S'$ indexed by $w_0$-conjugacy orbits on $S$.
At the level of Artin groups,
one shows (see for example \cite{jmichel}) that
$$A(W',S') \simeq C_{A(W,S)}(\Delta),$$
which is an algebraic reformulation of the case $p=1$, $q=2$ of
Question 3 (applied to $W\backslash V^{\reg}$, as in Example \ref{myexample}).

\begin{example}
\label{example04}
When $A(W,S)$ is a Artin group of type $E_6$, the $\Delta$-conjugacy
action is the non-trivial diagram automorphism and the centraliser
is an Artin group of type $F_4$.
$$
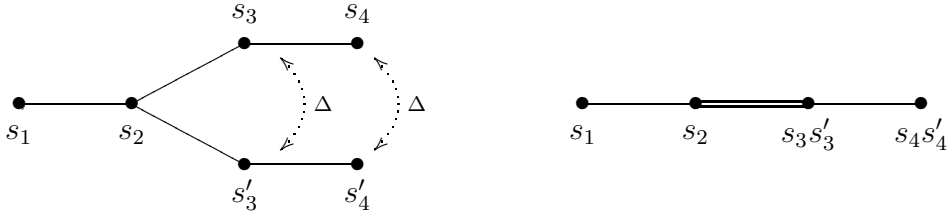
\xy 
(0,0)*+{\bullet}, (15,0)*+{\bullet},
(30,8)*+{\bullet}="c", (45,8)*+{\bullet}="a",
(30,-8)*+{\bullet}="d", (45,-8)*+{\bullet}="b",
(0,-4)*+{s_1}, (15,-4)*+{s_2},
(30,12)*+{s_3}, (45,12)*+{s_4},
(30,-12)*+{s_3'}, (45,-12)*+{s_4'},
\ar @{-} (0,0);(15,0),
\ar @{-} (15,0);(30,8),
\ar @{-} (30,8);(45,8),
\ar @{-} (15,0);(30,-8),
\ar @{-} (45,-8);(30,-8),
\POS"a" \ar @{<.>} @/^3ex/ "b"^\Delta,
\POS"c" \ar @{<.>} @/^3ex/ "d"^\Delta
\endxy
\qquad \qquad
\xy 
(0,0)*+{\bullet}, (15,0)*+{\bullet},
(30,0)*+{\bullet}="c", (45,0)*+{\bullet}="a",
(0,-4)*+{s_1}, (15,-4)*+{s_2},
(30,-4)*+{s_3s_3'}, (45,-4)*+{s_4s_4'},
\ar @{-} (0,0);(15,0),
\ar @{=} (15,0);(30,0),
\ar @{-} (30,0);(45,0)
\endxy$$
\end{example}
The main strategy throughout this article is to construct Garside 
structures with sufficient symmetries, so that centralisers of periodic
elements can be
computed as easily as in Example \ref{example04}.

Birman-Ko-Lee showed that the classical braid group $B_n$ admits, in addition
to the type $A_{n-1}$ Artin group structure, another Garside group structure
where the Garside element is a rotation $\delta$ of angle $2\pi\frac{1}{n}$.
In \cite{bdm}, we used this Garside structure to compute the centralisers
of powers of $\delta$, which solves Question 3 for $B_n$ and $q|n$.
Thanks to some rather miraculous diagram chasing, we were also able to
obtain the remaining case $q|n-1$.

Whenever $G$ is a group and $\Delta\in G$ is the Garside element of a 
certain Garside structure, the
centraliser $C_{G}(\Delta)$ is again a Garside group, and is easy to compute.
Note that the notion of \emph{periodic element} may be extended to this
setting: we say that $\rho\in G$ is \emph{$\frac{p}{q}$-periodic} if
$$\rho^q=\Delta^p.$$
This of course is relative to the choice of a particular Garside structure.
However,
for $B_n$, the Artin Garside element $\Delta$ and the Birman-Ko-Lee Garside
element $\delta$  are \emph{commensurable}
$$\Delta^2=\delta^n.$$
In particular, whether a given
 element $\rho$ is periodic or not does not depend on the choice between
the two standard Garside structures
(although the actual $p$ and $q$ may vary).

As the above particular case illustrates, it is very easy to compute
centraliser of Garside elements in Garside groups. In particular,
when trying to answer Questions 3 and 3' in a space whose fundamental
group is a Garside group, it is very tempting to expect to build a proof
on a positive answer to:

\begin{quote}
Question 4. \emph{Let $G$ be a Garside group with Garside element $\Delta$.
Let $\rho\in G$ be a periodic element with respect to $\Delta$.
Does $G$ admit a Garside structure with Garside element $\rho$?
}
\end{quote}

Note that a
positive answer to Question 4 would imply that the following question
also admits a positive answer:

\begin{quote}
Question 5. \emph{Let $G$ be a Garside group with Garside element $\Delta$.
Let $\rho\in G$ be a
periodic element with respect to $\Delta$.
Is the centraliser
$C_G(g)$ a Garside group?}
\end{quote}

Birman-Ko-Lee's discovery left many people puzzled: why should
there be \emph{two} natural Garside structures on $B_n$, how many
more Garside structures on $B_n$ remain to be discovered? Question 4
asks for a natural explanation in terms of periodic elements.
Although it was never written down as an official conjecture, there
was some initial hope that the answer might be positive. When $(W,S)$
is a finite Coxeter system, $A(W,S)$ admits
a \emph{dual Garside structure}, \cite{dualmonoid}, generalising
Birman-Ko-Lee's construction and giving partial answers when $q$ divides
$h$, the Coxeter number of $W$.

Our answer to Questions 4 and 5 is ``\emph{Almost!}''. Although we do
not have definite counterexamples, we do not expect them
to have positive answers, because we do not think that they are
phrased in a natural language. As the current article will illustrate,
the notion of
\emph{Garside groups} is artificially restrictive and
unstable under several basic operations. One should rather work
with \emph{Garside groupoids} and
\emph{weak Garside groups}. This generalisation was recently introduced
by Krammer and a related notion was independently studied by Digne-Michel,
\cite{krammer,dm}.
Just like Garside groups are groups of fractions of \emph{Garside monoids},
Garside groupoids are obtained by localising \emph{Garside categories}.
A monoid is a category with a single object, and a Garside monoid is a
Garside category with a single
object. Rewriting the whole theory of Garside groups
into categorical language is a surprisingly pleasant translation exercise:
everything works fine at essentially no cost. The additional syntactic
constraints are actually helpful, as they provide a quick test for
general statements about Garside groups: if a statement
cannot be ``categorified'', then it is probably false.

A \emph{categorical Garside structure} is a triple
$(\CC,\phi,\Delta)$ where $\CC$ is a small category,
$\phi$ is an automorphism of $\CC$ (it replaces the conjugacy action
of the Garside element, and should be thought of as a \emph{diagram
automorphism}) and $\Delta$ is a natural transformation from
the identity functor to $\phi$ (the family of \emph{Garside elements}), subject
to certain axioms. 
A \emph{Garside category} is a category $\CC$ which is part of such a triple.
It is cancellative and embeds in its Garside groupoid, the category $\CG$
obtained by adding formal inverses to all morphisms. A \emph{weak Garside
group} is the automorphism group $\CG_x$ of some object $x$ of a Garside
groupoid. When the category has a single object, one recovers the usual
notion of Garside group.

A Garside groupoid is \emph{cyclic} if the automorphism $\phi$ has
finite order.
Our ``almost answers'' to Questions 4 and 5 are as follows
(for a more precise phrasing of Theorem \ref{theo06}, see Theorems
\ref{theo94} and 
\ref{theo101}):

\begin{theo}
\label{theo06}
Let $\CG$ be a cyclic Garside groupoid, let $\gamma$ be a periodic loop
in $\CG$. Then
there exists a Garside groupoid $\CG_q$, together with an equivalence
of categories $$\Theta_q:\CG \to \CG_q$$
such that $\Theta_q(\gamma)$ is conjugate to a Garside element.
\end{theo}

\begin{coro}
The centraliser of a periodic element in a weak Garside group
is a weak Garside group.
\end{coro}

Note that, even when $\CG$ is a group, the groupoid $\CG_q$ may
have several objects. For the issues discussed here, it
is unavoidable to think in terms of groupoids.

These results have a geometric interpretation.
One associates to any Garside groupoid a simplicial classifying space
$\CN_{\Delta}\CG$,
which we call the \emph{Garside nerve}. When $\phi$ is trivial,
one shows that the Garside nerve is a cyclic set, in the sense of Connes.
More generally, when the Garside groupoid is cyclic, the Garside nerve
is very close to being a cyclic set: it is a $\Lambda_k^{\op}$-object
in the
category of sets, in the sense of B\"okstedt-Hsiang-Madsen. As a consequence,
the realisation $|\CN_\Delta\CG|$ comes equipped with a natural $S^1$-structure.
Theorem \ref{mykerekjarto} should be understood in that setting.

Our main construction, in Section \ref{section9},
is a sort of \emph{barycentric subdivision} operation for
Garside categories. 
For any Garside category $\CC$ with groupoid $\CG$,
we construct \emph{$m$-divided} Garside categories $\CC_m$ and
groupoids $\CG_m$, one for each integer $m\geq 1$.
We have $\CC_1=\CC$, but higher values of $m$ give
categories
$\CC_m$ with more objects than $\CC$. In particular, $m$-divided
categories of Garside monoids are usually not monoids.
The Garside nerve of $\CG_m$ is obtained
from that of $\CG$ by applying B\"okstedt-Hsiang-Madsen
$m$-subdivision functor, \cite{bhm}, and
their geometric realisations are homeomorphic.
This implies that there is
an equivalence of categories between $\CG_m$ and $\CG$.

In other words, when working with weak Garside groups, one may replace
$\CG$ by $\CG_m$. What is gained in the procedure is that the automorphism
group of $\CG_m$ is larger than that of $\CG$: the diagram automorphism
of $\CG_m$ can be thought of as an $m$-th root of that of $\CG$.
This, together with earlier ideas of Bestvina, is the main ingredient
in the proof of Theorem \ref{theo06}.

{\bf \flushleft Applications.}
As mentioned above (Example \ref{myexample}), we use our techniques
in a separate article, \cite{kapiun}, to obtain a complete description
of periodic elements, their conjugacy classes and their centralisers,
in braid groups associated with well-generated complex reflection groups.
Even in the case of spherical type Artin groups, these results are
new. However, the main application so far is
as an ingredient in the proof of $K(\pi,1)$ conjecture for complex reflection 
arrangements, which is the main result in \cite{kapiun}.
It is also likely that our results have algorithmic applications
to the conjugacy problem in Garside groups, in the spirit of
the work of Birman, Gebhardt and Gonz\'alez-Meneses, \cite{bgg}.

{\bf \flushleft Structure of this article.}
The first four sections provide basic terminology about Garside
categories, Garside groupoids and Garside germs.
They are included for the convenience
of the reader, to serve as a travel guide to a newborn theory,
rather than as an encyclopaedia.
Proofs are omitted and some axioms are stricter than actually needed.
A unified toolbox for Garside categories is currently being developed, in
collaboration with Fran\c cois Digne, Daan Krammer
and Jean Michel, and should eventually provide complete reference for the
missing details.
The next four sections cover more advanced material on Garside groupoids:
conjugacy problem (Section \ref{section5}), Galois theory (Section
\ref{section6}), classifying spaces (Section \ref{section7}) and
(a very brief account of a corollary of) Bestvina's approach
to non-positively curved aspects (Section \ref{section8}). Here
again, most proofs are omitted, because writing the details would probably
take over 30 additional pages and merely consist of straightforward
rephrasings of standard pieces from the theory of Garside groups,
copy-pasted from the works of Garside, ElRifai-Morton, Picantin,
Franco and Gonz\'alez-Meneses, Bestvina, Charney-Meier-Whittlesey and others.
The real core of this article consists of
Sections \ref{section9}, \ref{section10}, \ref{section11} --
this is where genuinely new material is introduced, that does not 
resemble anything that can be done without the categorical viewpoint.
The last two sections are devoted to easy illustrative examples.

\section{Graphs, germs and free categories}

Let $S$ be an oriented graph.
For any two vertices $x$ and $y$, denote by $S_{x\to}$ the set
of edges with source $x$, by $S_{\to y}$ the set of edges of
target $y$, and set $S_{x\to y}:=S_{x\to}\cap S_{\to y}$.

We may think of a (small) category $\CC$ as a graph $C$, whose vertices
are $\CC$-objects and edges are $\CC$-morphisms, together with a composition
law. In that respect, our notations
$$C_{x\to}, \; C_{\to y}, \; C_{x\to y}$$ are synonyms for 
$$\Hom_{\CC}(x,-), \; \Hom_{\CC}(-,y),\; \Hom_{\CC}(x,y).$$
Both notation systems will be used in the sequel.

The notion of \emph{germ} generalises the notion of \emph{category}
by allowing the composition law to be only partially defined:
\begin{defi}
A \emph{germ of (small) category} (or simply a \emph{germ}) is a pair
$\CS=(S,m)$ where $S$ is an oriented graph
together with partial ``composition'' maps
 $$m_{x,y,z}:S_{x\to y} \times S_{y\to z} \to S_{x\to z},$$
one for each triple $(x,y,z)$ of vertices,
satisfying the following axioms:
\begin{itemize}
\item[(assoc)] let $x,y,z,t$ be vertices; the two natural partial maps
$$m_{x,y,t}\circ (1_{S_{x\to y}}\times m_{y,z,t}):
S_{x\to y} \times S_{y\to z} \times S_{z\to t} \to S_{x\to z}$$ and
$$m_{x,z,t}\circ (m_{x,y,z} \times 1_{S_{z\to t}}):
S_{x\to y} \times S_{y\to z} \times S_{z\to t} \to S_{x\to z}$$
coincide (in particular, we ask for these maps to have the same domain
of definition).
\item[(unit)] for all vertex $x$, there exists an element $1_x\in S_{x\to x}$
such that, for all vertices $x,y$, the partial maps
$$m_{x,x,y}(1_x,\cdot):S_{x\to y} \to S_{x\to y}$$
and
$$m_{x,y,y}(\cdot,1_y):S_{x\to y} \to S_{x\to y}$$
coincide with the identity map $1_{S_{x\to y}}$
(in particular, we ask for these maps to be everywhere defined).
\end{itemize}
\end{defi}

It is an easy check that (assoc) implies that if $s_1,\dots,s_n\in S$,
the choice of a complete bracketing of $s_1\dots s_n$
has no impact in whether the product is defined or not, nor on the
occasional value of this product. In (unit), the element $1_x$ is unique.

\begin{example}
A germ where $m_{x,y,z}$ is everywhere defined is nothing
but a category.
\end{example}

Let $\CS:=(S,m)$ be a germ.
Let $S^{*}$ be the category of walks on $S$:
its objets are the vertices of $S$, and a morphism from $x$ to $y$ is
a finite sequence (possibly empty) of $S$-edges $(s_1,\dots,s_k)$ such
that $s_1\in S_{x\to }$, $s_k\in S_{\to y}$, and the source of $s_{i+1}$
is the target of $s_i$. Paths are composed by concatenation.
Say that two paths $(s_1,\dots,s_k)$ and $(t_1,\dots,t_l)$ are
\emph{elementarily $\CS$-equivalent}, which is denoted by
$$(s_1,\dots,s_k)\underset{1}{\sim}(t_1,\dots,t_l),$$
if one (or both) of the following conditions is satisfied
\begin{itemize}
\item[(I)]There exists $i$ such that
 $(s_i,s_{i+1})$ is in the domain of definition of $m$ and
$$(t_1,\dots,t_l) = (s_1,\dots,m(s_i,s_{i+1}),\dots,s_k).$$
\item[(II)] There exists $i$ such that $s_i$ is the unit of some object $x_i$
and 
$$(t_1,\dots,t_l) = (s_1,\dots,s_{i-1},s_{i+1},\dots,s_k).$$
\end{itemize}

\begin{remark}
\begin{itemize}
\item[(i)] In both cases, we have $l=k-1$.
The second case is \emph{almost} a particular
case of the first one, except that it allows
 $() \underset{1}{\sim} (1_{x})$ which is not covered by the first case.
\item[(ii)] In point (i) above, the notation $()$ is slightly ambiguous,
since it refers to the trivial walk \emph{starting at $x$}. When needed,
we will use the notation $\phantom{.}_x()$ to lift ambiguity.
\item[(iii)]  Section 0 of
my earlier paper \cite{dualmonoid} gives
 a monoid version of $\underset{1}{\sim}$ but contains a mistake,
pointed out by Deligne: case (II) was forgotten, which had the
unpleasant effect of adding artificial units.
\end{itemize}
\end{remark}

Let $\sim$ be the reflexive symmetric transitive closure
of $\underset{1}{\sim}$. It is clear that the concatenation of paths is
compatible with $\sim$. We obtain a category
$S^*/\sim$ whose morphisms are $\sim$-equivalence classes of paths.

\begin{defi}
The \emph{free category on a germ $\CS$}, denoted by $\CC(\CS)$, is the
quotient category $S^*/\sim$.
\end{defi}

The terminology is justified by the fact that $\CS\mapsto \CC(\CS)$
is a left adjoint of the
forgetful functor from the category of small categories to the
category of germs.

\begin{notation}
When there is no ambiguity, we simply denote a germ $\CS=(S,m)$ by
its underlying graph $S$. Furthermore, we use the notation $s\in S$
to mean that \emph{$s$ is an oriented edge of $S$} or, with equivalent
categorical language, that \emph{$s$ is a morphism of $S$}.
We say that $s\in S$ is an \emph{object} if it is the identity morphism
$1_x$ associated with a vertex $x$ (as with categories, one may think of
objects as particular morphisms). If $s,t,u\in S$, the notation
$st=u$ means that ``$(s,t)$ is in the definition domain of $m$
and $m(s,t)=u$.''
\end{notation}

The partial product endows $S_{x\to}$ with a poset structure:
$$s'\leq s  \stackrel{\text{def}}{\Longleftrightarrow} \exists s'',s's''=s.$$
Similarly, it endows $S_{\to y}$ with a poset structure:
$$s \geq s''  \stackrel{\text{def}}{\Longleftrightarrow} \exists s',s's''=s.$$

\section{Garside categories}
\label{sectiongarside}

A monoid is a category with a single object.
While Garside monoids have been studied for quite a while, Garside
categories are a recent invention. Krammer gives a very neat definition
in \cite{krammer}. As he points out, many theorems about Garside
monoids, as well as their proofs, may be rewritten with no effort
in the context of Garside categories. In fact, the categorical viewpoint is
arguably simpler and more natural, even when dealing with Garside monoids.
When Dehornoy was illustrating his articles and lectures
with pictures where elements of Garside monoids were represented by arrows,
he was adopting a categorical viewpoint without realising
it\footnote{or without saying it...}.
The categorical viewpoint is also implicit in Deligne's article
\cite{deligne} (see Example \ref{exampledeligne} below).

As with Garside monoids, one may define a Garside category by generators
and relations (as Krammer does), or by showing that it satisfies
certain axioms. Also, it is tempting to relax some of the axioms, to 
allow for ``quasi''-Garside categories (where ``quasi'' could mean several
different things) retaining some of the main properties. It not
yet clear which final optimal axiom set will be retained, and what
follows is a simple variant.

The starting kit is a \emph{basic triple}
$$(\CC,\CC\stackrel{\phi}{\to}\CC,1_{\CC}\stackrel{\Delta}{\Rightarrow}\phi)$$
where $\CC$ is a small category, $\phi$
is an automorphism\footnote{It is possible to do part of the theory
by simply assuming that $\phi$ is an endofunctor.} 
of $\CC$ and $\Delta$ is a natural
transformation from the identity functor to $\phi$.

\begin{example}
Let $M$ a Garside monoid with Garside element $\delta$, 
view it as a category with one point $\ast$ and
arrows labelled by elements of $M$. For any arrow $m\in M$,
take $\phi(m):=\delta^{-1}m\delta$, and take $\Delta$ to be the right
multiplication by $\delta$. The naturality of $\Delta$ is expressed
by $\forall m\in M, \delta \phi(m) = m \delta$.
\end{example}

To be consistent with Dehornoy's pictures and other classical material,
we use some conventions.
Arrows in $\CC$ compose like paths in algebraic topology:
$x\stackrel{f}{\to}y \stackrel{g}{\to}z$ is composed into
$x\stackrel{fg}{\to}g$.
The functor $\phi$ is  denoted as if it was a ``right
conjugacy action'':
the image of $x\stackrel{f}{\to}y$ by $\phi$ is
$x^{\phi} \stackrel{f^{\phi}}{\to}y^{\phi}$.
When $x\in\CC$, $\Delta$ gives a morphism $x\stackrel{\Delta(x)}{\to} x^{\phi}$
which we like to simply denote by $\Delta_x$ or even $\Delta$, calling it
``the $\Delta$ of $x$''. Any element $\Delta_x$ should be called
a \emph{Garside element}.

Krammer asks for $\CC$ to have a finite
number of objects and a finite number of atoms.
If one wants to study infinite type Artin groups, infinite number of atoms
should be allowed. And since our main construction turns atoms into objects
of newer categories, it is natural to allow infinite number of objects.
Krammer's
axiom (GA2) exactly expresses that $\Delta$ is a natural transformation.

\begin{defi}
\label{defisimple}
Let $(\CC, \phi,\Delta)$ be a basic triple.
A morphism $x \stackrel{f}{\to} y$ is
\emph{simple} if there exists a morphism
$y \stackrel{\overline{f}}{\to} x^{\phi}$ such that $f\overline{f}=\Delta_x$.
We denote by $S(\CC, \phi,\Delta)$ (or simply $S$) the \emph{graph
of simples}, the subgraph of $\CC$ whose edges are simple morphisms.
\end{defi}

A category is \emph{connected} if the underlying graph in connected
(in the unoriented sense).

A category is \emph{atomic} if for any morphism $f$, there is a bound
on the length $n$ of a factorisation $f=f_1\dots f_n$, where the
$f_i$ are non-identity morphisms.
This implies that there are no non-trivial invertible morphisms;
in particular, whenever a limit or colimit exists, it
is \emph{unique} (really unique, not just up to automorphism).
A nontrivial morphism $f$ which cannot be factorised into two
nontrivial factors is an \emph{atom}.

A category is \emph{(weighted) homogeneous} if there exists a length
function $l$ from the set of $\CC$-morphisms to $\BZ_{\geq 0}$
such that $l(fg)=l(f)+l(g)$ and
$(l(f)=0)\Leftrightarrow (\text{$f$ is a unit})$. It is clear that
homogeneous categories are atomic.

\begin{defi}
If $\CC$ is a category, we denote by $A(\CC)$ (or simply $A$)
the \emph{atom graph} of $\CC$, the subgraph of $\CC$ whose edges
are atoms.
\end{defi}
In an atomic category, any morphism is a product of
atoms; in other words, $A(\CC)$ generates $\CC$.

A category is \emph{cancellative} if,
whenever a relation $afb=agb$ holds between composed morphisms,
it implies $f=g$. Very often, this implies that there is at most one way
to add an certain arrow to a commutative diagram.
In the context of Definition \ref{defisimple}, this implies that
$\overline{f}$ is unique.

\begin{defi}
A \emph{(homogeneous)
categorical Garside structure} is a triple $(\CC,\phi,\Delta)$
such that:
\begin{itemize}
\item $\CC$ is a category, $\phi$ an automorphism of $\CC$ and 
$\Delta$ a natural transformation from the identity functor to $\phi$,
\item $\CC$ is homogeneous
\footnote{\label{foot}This condition is too restrictive, but sufficient for
most applications, in particular when considering Garside categories
naturally associated
with braid groups, like in \cite{kapiun}. As in \cite{dualmonoid}
and in \cite{krammer}, one can replace homogeneity by atomicity. This
may yet not be the ultimate general condition. To avoid technicalities,
we stay with this overly conservative axiom.} and cancellative,
\item atoms are simple: $A(\CC)\subseteq S(\CC,\phi,\Delta)$,
\item for all object $x$, $(\CC_{x \to},\leq)$ and $(\CC_{\to x},\geq)$
are lattices.
\end{itemize}
It has \emph{finite type} if $S(\CC,\phi,\Delta)$
is finite.
Let $k\in \BZ_{\geq 1}$. We say that  $(\CC,\phi,\Delta)$  is
\emph{$k$-cyclic} if $\phi^k=1$. It is 
\emph{cyclic} if it is $k$-cyclic for some $k\geq 1$. 
\end{defi}

Note that we do not require Garside categories to be connected nor
to be non-empty.

\begin{defi}
A \emph{Garside category} is a category $\CC$ that may be equipped
with $\phi$ and $\Delta$ to obtain a categorical Garside structure.
\end{defi}

Krammer's definition of Garside categories corresponds to our
\emph{finite type} Garside categories. \emph{Infinite type} Garside
categories fail to satisfy property (P7) from Krammer's Theorem 36
(automaticity), but satisfy the remaining properties.

\begin{example}
The basic triple associated to a Garside monoid is a finite
type Garside triple. When considering ``quasi-Garside'' monoids
(with infinite number of simples, such as in \cite{digne} or \cite{free}),
one obtains an infinite type Garside triple.
This also provides us with examples of non-cyclic Garside categories.
\end{example}

\begin{defi}
A Garside groupoid is a groupoid obtained by adding formal inverses
to all morphisms in $\CC$, where $\CC$ is a Garside category.
\end{defi}

\begin{defi}
Let $(\CC,\phi,\Delta)$ be a Garside category, with Garside
groupoid $\CG$.
The \emph{structure group of $\CG$} at an object $x$ is 
$$\CG_x:=\End_{\CG}(x,x).$$

A \emph{weak Garside group} is a group that is isomorphic to a
structure group of a Garside groupoid.
\end{defi}

The isomorphism type of $\CG_x$ only depends on the connected component
of $\CC$ containing $x$.

\begin{remark}
If $\CC$ has a single object (hence is a Garside \emph{monoid})
then $\CG_x$ is a Garside group, in the previously traditional sense.
Thus Garside groups are weak Garside groups.
Conversely, there is no good reason to expect all weak
Garside groups to be Garside groups. We give at the end of this
paper an example to illustrate this.
Because we think that \emph{weak Garside groups} are much more natural
to consider than traditional \emph{Garside groups}, we hope that
the terminology will evolve and that people will eventually call
\emph{Garside groups} what we call \emph{weak Garside groups}.
Until then, it is probably safer to keep the ``\emph{weak}.''
\end{remark}

\begin{defi}
Let $(\CC,\phi,\Delta)$ be a categorical Garside structure.

The \emph{Garside dimension} of $(\CC,\phi,\Delta)$ is the 
element of $\BZ_{\geq 0} \cup \{+\infty\}$ defined by
$$\dim_{\Delta}(\CC,\phi,\Delta) := \sup \{ n\in \BZ_{\geq 0} | \exists
\text{ simples } s_0,\dots,s_n \text{ such that } s_0 < s_2 < \dots < s_n \}.$$
\end{defi}

Let $(\CC,\phi,\Delta)$ be a categorical Garside structure
with associated groupoid $\CG$.
Let $f\in \CG$. The classical arguments from the theory of normal
forms in Garside groups are applicable here, and one sees
that there is a unique way to write $f$ as a product
$$f = s_1 s_2 \dots s_l \Delta^k,$$
where $s_1,\dots,s_l$ are simples with sources $x_1,\dots,x_l$, $k\in \BZ$,
and we have
$$\text{for all $i$, } s_i = s_is_{i+1} \wedge \Delta_{x_i}$$
and $$s_1 < \Delta_{x_1}.$$
In the above formulae, the symbol $\wedge$ refers to the left gcd,
\ie, the inf with respect to $\leq$.
Note that we choose to put $\Delta^k$ at the end rather than at
the beginning, to stay in line with
conventions used by Bestvina and others.

\begin{defi}
We say that $s_1 s_2 \dots s_l \Delta^k$ is the \emph{(left greedy) normal form}
of $f$. The integers $k$ and $k+l$ are respectively the 
\emph{infimum} and \emph{supremum} of $f$, denote by 
$\lfloor f \rfloor$ and $\lceil f \rceil$. The integer $l$ is the
\emph{canonical length} of $f$.
\end{defi}

When the Garside structure has finite type, this normal form
is part of an \emph{automatic structure} for $\CG$.

\begin{remark}
Let $k$ be a positive integer. Any $(\CC,\phi,\Delta)$ categorical Garside
structure on $\CC$ admits a ``$k$-th power'' $(\CC,\phi^k,\Delta^k)$, whose
 Garside elements are products of $k$ consecutive
$\Delta$'s. Our main construction below provides (a sort of) inverse
to this power operation on Garside structures.
\end{remark}

\section{Garside germs}

This section explains how to reconstruct a Garside category from a
\emph{Garside germ}, which should be thought of as
a tentative ``graph of simples''.

While Dehornoy-Krammer's syntactic approach is very effective at handling
abstract Garside categories whose graph of simples is hard to understand,
many situations provide us with a natural candidate for the graph of simples,
while there may be no natural complemented presentation and the cube
axiom may be unpractical to check.

For Garside monoids, a more intrinsic strategy is explained in my
joint work with Fran\c cois Digne and Jean Michel, \cite{bdm}. The
material presented
here generalises \cite{bdm} and is in line with the viewpoint of Digne
and Michel on \emph{locally Garside categories} (a weaker notion),
\cite{dm}.

Let $(\CC,\phi,\Delta)$ be a categorical Garside structure.
Let $A$ be the
graph of atoms and $S$ be the graph of simples of $\CC$.
The automorphism $\phi$ induces automorphisms of these graphs.

Krammer explains how to reconstruct $\CC$ as a quotient
of the path category of the atom graph, via the
choice, for each pair $a,b$ of atoms, of paths $\underline{a\backslash b}$
and $\underline{b\backslash a}$ (a path in the atom graph is a formal
sequence of atoms),
corresponding to atomic decompositions in $\CC$ of the
right factors $a\backslash b$ and
$b\backslash a$ of the colimit $a\vee b= a(a\backslash b)= a(b\backslash a)$.
He actually goes the other way around: he starts with an abstract graph,
together with choices of walks $\underline{a\backslash b}$ and
$\underline{b\backslash a}$, choices of walks expressing each $\Delta_x$,
and an automorphism of the whole structure;
when these choices satisfy certain conditions (Dehornoy's
cube condition, condition for naturality of $\Delta$,...),
then he declares that the quotient of the path category of the abstract graph
modulo the relations
$a (\underline{a\backslash b}) = b(\underline{b\backslash a})$ is a Garside
category.

We prefer to view $\CC$ as the free category on the germ 
$\CS:=(S,m)$, where $S$ is the graph of simple and $m$ is the
restriction of the $\CC$-composition law
$$S_{x\to y} \times S_{y\to z} \to \CC_{x\to z}$$
to the preimage of $S_{x\to z}$. We call $\CS$ the \emph{germ of simples}
of  $(\CC,\phi,\Delta)$.

\begin{lemma}
Let $\CS$ be the germ of simples of a categorical Garside
structure $(\CC,\phi,\Delta)$.
The natural morphism $\CC(\CS)\to \CC$ is an isomorphism.  
\end{lemma}

Conversely, suppose we are given a germ $\CS$. How can we know
whether it is the germ of simples of a Garside category?
The axioms of Garside categories may be rewritten into a set 
of axioms characterising such germs.

\begin{defi}
A germ $\CS=(S,m)$ is \emph{(homogeneous) Garside} if the following conditions
are satisfied:
\begin{itemize}
\item[(i)] It is homogeneous\footnote{See footnote
\ref{foot} on page \pageref{foot}.
Note also that the definitions
for atomicity, homogeneity and cancellativity easily generalise
from categories to germs (we leave the details to the reader).}
and $\CC$-cancellative\footnote{This notion is the natural generalisation
of the $\mathbf{M}$-cancellativity from \cite[Section 0]{dualmonoid}}.
\item[(ii)] For any vertex $x$, $(S_{x\to},\leq)$ admits a maximal element
$\Delta_x$.
\item[(iii)] Denote by $x^{\phi}$ the target of $\Delta_x$. The
map
\begin{eqnarray*}
(S_{x\to },\leq) & \to & (S_{\to x^{\phi}},\geq) \\
s &\mapsto & \overline{s} \text{ such that } s\overline{s}=\Delta_x
\end{eqnarray*}
is an isomorphism.
\item[(iv)] For all $x$, the poset $(S_{x\to},\leq)$ is a lattice.
\end{itemize}
\end{defi}

Note that, in (iii), the map $s\mapsto \overline{s}$ is well-defined
because $\CS$ is cancellative. Because of (iii),
$\Delta_x$ is also the maximal element of $(S_{\to x^{\phi}},\geq)$
and the apparent chirality of the axiom set is only an optical
illusion.

For any objects $x,y$, the
isomorphism $(S_{x\to },\leq) \simeq (S_{\to x^{\phi}},\geq)$
restricts to a bijection $S_{x\to y} \simeq S_{y\to x^{\phi}}$.
Applying this twice, we obtain a bijection
$$S_{x\to y} \simeq S_{y\to x^{\phi}} 
\simeq S_{x^{\phi} \to y^{\phi}}$$
which may be shown to be part of a global isomorphism $\phi:\CS\to \CS$.
It extends to an automorphism $$\phi:\CC(\CS)\to \CC(\CS)$$ such that
$$\Delta:x\mapsto (\Delta_x)/\sim$$ is a natural transformation
$1\stackrel{\Delta}{\Rightarrow}\phi$.

\begin{theo}
The germ of simples of a categorical Garside structure is
a Garside germ.
Conversely, if $\CS$ is a Garside germ, then $(\CC(\CS),\phi,\Delta)$,
where $\phi$ and $\Delta$ are as constructed above, is a categorical
Garside structure whose germ of simples is isomorphic to $\CS$.
\end{theo}

\begin{example}
\label{exampledeligne}
Let $\CA$ be a finite real reflection arrangement.
Let $ch(\CA)$ be the
set of chambers (connected components of the complement of the
reflecting hyperplanes). There is a natural distance on $ch(\CA)$ ($d(C,C')$
is the number
of walls separating the two chambers).  Let $\CS$ be the germ whose underlying
graph is the complete oriented graph on $ch(\CA)$ (edges are pairs
$(C,C')$ of chambers) and such that $(C,C')(C',C'')$ is defined
as equal to $(C,C'')$ when $C'$ lies on a geodesic from $C$ to $C''$ (and
not defined otherwise). The category  
$\CC(\CS)$ is isomorphic to the category denoted by
 $\Gal_+$ in \cite{deligne}.
Deligne shows that if $\CA$ is \emph{simplicial},
e.g. if it is the reflection arrangement of a finite real reflection
group, then $\CS$ is a Garside germ and $\CC(\CS)$ is a Garside
category.
\end{example}

\section{Automorphisms of Garside categories}

\begin{defi}
Let $(\CC,\phi,\Delta)$ be a categorical Garside structure.
An \emph{automorphism} of $(\CC,\phi,\Delta)$ is an
automorphism $\psi$ of $\CC$ such that
$$\phi\psi = \psi \phi$$
and, for all object $x$,
$$\psi(\Delta_x)=\Delta_{\psi x}.$$
\end{defi}

If $(\CC,\phi,\Delta)$ is a categorical Garside structure,
the Garside automorphism
$\phi$ is an automorphism of $(\CC,\phi,\Delta)$:
indeed, if $x$ is an object, the naturality
$1\stackrel{\Delta}{\Rightarrow} \phi$
applied to the morphism $\Delta_x$ gives the commutative diagram
$$\xymatrix{ x \ar[r]^{\Delta_x} \ar[d]_{\Delta_x} &
x^{\phi}\ar[d]^{(\Delta_x)^{\phi}} \\
x^{\phi} \ar[r]_{\Delta_{x^{\phi}}} & x^{\phi^2} }$$
and, by cancellativity, 
$$(\Delta_x)^{\phi} =\Delta_{x^{\phi}}.$$

\begin{theo}
Let $(\CC,\phi,\Delta)$ be a categorical Garside structure.
Let $\psi$ be an automorphism of $\CC$. Then
$(\CC^{\psi},\phi |_{\CC^{\psi}},\Delta|_{\CC^{\psi}})$ is a categorical
Garside structure. 
\end{theo}

By $\CC^{\psi}$, we mean the subcategory of $\CC$ consisting of
morphisms invariant under $\psi$.

\begin{proof}
Since $\CC^{\psi}$ is a subcategory of $\CC$, it is atomic and
cancellative.

Let $x\in\CC^{\psi}$ be an object.
Since ${\psi}(x)=x$, we have
$\psi(\Delta_x) =\Delta_{{\psi}(x)}=\Delta_x$.
This justifies that $\Delta|_{\CC^{\psi}}$ is indeed a natural 
transformation in $\CC^{\psi}$.

Let $a$ be an atom of $\CC$, let $x$ be the source of $a$.
Assume that $\psi(x)=x$.
For any integer $k\geq 0$, let $\psi_ka$ be the colimit
of $\{a,\psi(a),\psi^2(a),\dots, \psi^k(a)\}$.
By atomicity,
since $\psi_ka\leq \Delta_x$, the sequence $a,\psi_1a,\psi_2a,\dots$
is eventually constant. Its limit $\psi_*a$ is $\psi$-invariant.

Any atom $b$ of $\CC^{\psi}$ is obtained this way: if $a$ is a $\CC$-atom
such that $a\leq b$, then $\psi_*a\leq b$, thus $\psi_*a=b$. Since
$\psi_*a\leq \Delta_x$, atoms are simple.

If $f,g\in \CC^{\psi}$ have the same source, then 
they admit a $\CC$-colimit $f\vee g$. This colimit divides a certain power
of $\Delta$. As above, this implies that the infinite family
$$\{ f\vee g,\psi(f\vee g),
\psi^2(f\vee g), \dots,\}$$ admits a colimit, which clearly is a
$\CC^{\psi}$-colimit for $f$ and $g$.
\end{proof}

\begin{coro}
\label{corocommute}
Let $(\CC,\Delta,\phi)$ be a categorical Garside structure.
Let $x$ be an object of $\CC$.
Let $p\in \BZ_{\geq 0}$ and consider the element $\Delta^p$
with source $x$. Assume that this element is a loop, namely
that $\Delta^p \in \CG_x$. Then the centraliser 
$C_{\CG_x}(\Delta^p)$ is a weak Garside group, namely the
structure group at $x$ of the Garside structure
$(\CC^{\phi^p},\phi\vert_{\CC^{\phi^p}},\Delta\vert_{\CC^{\phi^p}})$.
\end{coro}

\begin{proof}
Because $1\stackrel{\Delta}{\Rightarrow} \phi$, we have
$c\Delta^p = \Delta^p c^{\phi^{p}}$.
Thus,
for any $c\in \CG_x$, the conditions $c\in C_{\CG_x}(\Delta^p)$ 
and $c\in \CG^{\phi^p}$ are equivalent.
\end{proof}

\begin{remark}
Even when $\CC$ is connected, the category $\CC^{\psi}$ may be disconnected.
\end{remark}

\section{Loops and summits}
\label{section5}

Let $\CG$ be a category.
Let $g,g',c\in\CG$.
We write $$g'=g^c$$ as a synonym for the relation $$gc=cg'.$$
By syntactic constraints, this cannot happen unless $g$ and $g'$ are
\emph{loops}, \ie, if $g\in\CG_{x\to x}$ and $g'\in \CG_{y\to y}$
for some objects $x,y$.

\begin{defi}
We say that two loops $g$ and $g'$ are \emph{conjugate}, and denote this
by $g\sim g'$ if there exists $c$ such that $g^c=g'$.
\end{defi}

If $\CG$ is a groupoid, then $\sim$ is an
equivalence relation.

\begin{defi}
Let $\CG$ be a groupoid. The \emph{conjugacy category}
of $\CG$ is the category $\Omega\CG$ whose object set is
$$\bigsqcup_{x \text{ $\CG$-object}} \CG_x$$
and such that $\Hom_{\Omega\CG}(g,g'):=\{c\in \CG| g^c=g'\}$.
\end{defi}

The composition law is the obvious one:
$$(g\stackrel{c}{\to}g')\cdot 
(g'\stackrel{c'}{\to}g''):=g\stackrel{cc'}{\to}g''.$$
The conjugacy classes of $\CG$ are the connected components of $\Omega\CG$.

The conjugacy category is clearly a groupoid.
If $g\in \CG_x$, we have $$C_{\CG_x}(g) \simeq (\Omega \CG)_g,$$
where $C_{\CG_x}(g)$ is the centraliser $\{c\in \CG_x | gc=cg\}$ and
$(\Omega \CG)_g$ is the structure group of $\Omega\CG$ at $g$.

Assume now that $\CG$ is the Garside groupoid of a
categorical Garside structure $(\CC,\phi,\Delta)$.
Let $g\in \CG$. Write $g$ in normal form
$$g=s_1\dots s_l \Delta^k.$$
Recall that the
\emph{infimum} and \emph{supremum} of $g$ are,
respectively, the integers $\lfloor g \rfloor:=k$ and $\lceil g\rceil :=k+l$.

\begin{defi}
A loop $g$ in $\CG$ is a \emph{summit} if, 
for all $h$ such that $g\sim h$,
$$\lfloor h \rfloor \leq  \lfloor g \rfloor \quad \text{and} \quad
\lceil g \rceil \leq \lceil h \rceil.$$

The \emph{summit category} of $\CG$ is the full category
$\Omega_0\CG$ of $\Omega \CG$
whose objets are the summits.
\end{defi}

Any loop is conjugate to a summit. For any loop $g$, a certain procedure
called ``cycling and decycling'' (after \cite{em,picantin})
 may be applied to
obtain a $c$ such that $g^c$ is a summit (we won't use this procedure;
it easily generalises from the earlier settings).
The problem of determining
whether $g$ and $h$ are conjugate, or computing the centraliser of $g$,
can be reduced to the same problem dealing with summits rather than just loops.

The following key lemma traces back to Garside:

\begin{lemma}
\label{lemmasss}
Let $g \stackrel{c}{\to} h$ be a morphism in $\Omega_0\CG$.
Write $c=c_1\dots c_l\Delta^k$ in normal form.
Set $g_0:=g,g_1:=g^{c_1},g_2:=g^{c_1c_2},\dots,g_l:=g^{c_1\dots c_l}$.
Then $g_0,\dots, g_l$ are summits.
\end{lemma}

\begin{coro}
Let $\CG$ be a finite type Garside groupoid.
The conjugacy problem in $\CG$,
and the problem of finding presentations for the centraliser of a given loop,
are solvable.
\end{coro}

\section{Coverings}
\label{section6}

\begin{defi}
Let $\CG$ be a connected
groupoid, let $x_0$ be an object (``a basepoint'') of $\CG$. 
The \emph{universal cover} (aka \emph{Cayley category})
of $\CG$ with respect to $x_0$ is the category
$\widetilde{\CG}_{x_0}$ as follows:
\begin{itemize}
\item objects are $\CG$-morphisms $x_0\stackrel{f}{\to} x$ with source $x_0$;
\item for any pair
$x_0\stackrel{f}{\to} x$, $x_0\stackrel{g}{\to} y$ of objects, there exists
a unique $\widetilde{\CG}_{x_0}$-morphism from $f$ to $g$;
we denote this morphism by the formal symbol $f^{-1}g$.
\end{itemize}
\end{defi}

Morphisms compose the way they should: since
$\Hom_{\widetilde{\CG}_{x_0}}(f,h)$ contains a single
element, we have to set:
$$f^{-1} g\cdot g^{-1}h=f^{-1}h.$$
In particular, there is a ``covering'' functor
\begin{eqnarray*}
p:\widetilde{\CG}_{x_0} & \longrightarrow & \CG \\
x_0\stackrel{f}{\to}x & \longmapsto & x \\
f^{-1} g & \longmapsto & f^{-1} g
\end{eqnarray*}
Of course, $f^{-1} g$ is a formal symbol when viewed in the universal
cover, and an actual morphism $x\to y$ when viewed in $\CG$. 

The universal cover comes equipped with a natural basepoint $1_{x_0}$.
It is clear that applying again the universal cover construction
to  $\widetilde{\CG}_{x_0}$ does not bring anything new.

Assume now that $\CG$ is the Garside groupoid of a Garside category
$(\CC,\phi,\Delta)$. Let $x_0$ be an object of $\CC$.
Let $\widetilde{\CC}$ be the subcategory of $\widetilde{\CG}_{x_0}$
consisting of morphism whose image under $p$ lies in $\CC$.
For any object $x_0\stackrel{f}{\to} x$ of $\widetilde{\CC}$,
set $$f^{\widetilde{\phi}}:=f\Delta_x.$$
If $f\stackrel{f^{-1}g}{\to} g$ is a $\widetilde{\CC}$-morphism,
we set $$(f^{-1}g)^{\widetilde{\phi}}:=(f^{\widetilde{\phi}})^{-1}g
^{\widetilde{\phi}}.$$
One checks that $(f^{-1}g)^{\widetilde{\phi}}$ is indeed a
$\widetilde{\CC}$-morphism, that $\widetilde{\phi}$ is an automorphism
of $\widetilde{\CC}$, 
and that
$$\widetilde{\Delta}: (x_0\stackrel{f}{\to}x) \longmapsto (f\to f\Delta_x)$$
is a natural transformation from the identity functor to $\widetilde{\phi}$.

\begin{lemma}
The triple $(\widetilde{\CC},\widetilde{\phi},\widetilde{\Delta})$
is a Garside category, whose Garside groupoid is
$\widetilde{\CG}_{x_0}$.
\end{lemma}

More generally, we may construct an intermediate cover for each prescribed
isotropy subgroup of $\CG_{x_0}$.

\section{The Garside nerve}
\label{section7}

Let us begin by recalling the categorical viewpoint on group cohomology.
A good introduction may be found in \cite{segal}. 

The \emph{nerve} of a small category $\CC$ is the simplicial
set $\CN\CC$ whose $0$-skeleton is the object set of $\CC$
and whose $n$-simplices, $n\geq 1$, are sequences
$(f_1,\dots,f_n)$ of
$\CC$-morphisms composable as follows:
$$\xymatrix{
x_0 \ar[r]^{f_1} & x_1 \ar[r]^{f_2} & x_2 \ar[r]^{f_3}  & x_3  \ar@{..>}[r]
& x_{n-1} \ar[r]^{f_{n}} & x_{n}}.$$
Face maps correspond to removing objects (and composing or dropping
morphisms
accordingly) and degeneracy
maps correspond to inserting identity morphisms at a given object.

\begin{lemma}
If $\CG$ is a groupoid,
the geometric realisation $|\CN\widetilde{\CG}_{x_0}|$ is contractible.
\end{lemma}

\begin{proof}
This is because
the category $\widetilde{\CG}_{x_0}$ is equivalent to the trivial category
with only one arrow, and the functor sending a small category $\CC$
to the realisation of its nerve is a functor of $2$-categories: it sends
small categories to topological spaces, functors to continuous maps and
natural transformations of functors to homotopies between continuous maps.
\end{proof}

When $\CG$ has a single object $x_0$, \ie, when it is a group,
the $\BZ$-basis of the $\bar$ resolution is
indexed by $\CN\widetilde{\CG}_{x_0}$. 

In the general situation,
$\CN\widetilde{\CG}_{x_0}$ is related to $\CN\CG$ via the
``$\bar$ notation'':
let $x_0 \stackrel{g_0}{\to}x_1$ be an object of $\widetilde{\CG}_{x_0}$.
Let $g:=(g_1,\dots,g_n)$ be a $n$-simplex in $\CN\CG$ such that
the source of $g_1$ is $x_1$:
$$\xymatrix{
x_1 \ar[r]^{g_1} & x_2 \ar[r]^{g_2} & x_3 \ar[r]^{g_3}  & x_4  \ar@{..>}[r]
& x_{n} \ar[r]^{g_{n}} & x_{n+1}}$$
of $\CG$-morphism.
There exists a unique $n$-simplex $(h_1,\dots,h_n)$
in $\CN\widetilde{\CG}_{x_0}$
that satisfies the following conditions:
\begin{itemize}
\item $h_i$ is mapped to $g_i$ by the covering functor,
\item the source of $h_1$ is $g_0$.
\end{itemize}

Indeed, the only solution is that each $h_i$ should be the unique
$\widetilde{\CG}_{x_0}$-morphism from $g_0\dots g_{i-1}$ to
$g_0\dots g_{i}$.

Instead of using the ``morphisms'' notation $$(h_1,\dots,h_n),$$ we may
represent
this $n$-simplex by its $\bar$ notation $$g_0[g_1|\dots | g_n]$$
or its ``endpoints'' notation $$(g_0,g_0g_1,\dots,g_0g_1\dots g_n).$$

For $g\in \CG_{x_0}$ and $g_0[g_1|\dots | g_n] \in \CN\widetilde{\CG}_{x_0}$,
we set
$$g\cdot g_0[g_1|\dots | g_n]:=gg_0[g_1|\dots | g_n].$$
This defines a left action of the structure group $\CG_{x_0}$
on $\CN\widetilde{\CG}_{x_0}$. This action is free.
It is natural to think that $\CN\widetilde{\CG}_{x_0}$ is a
``$\bar$ resolution for groupoids''.

We don't know if the following terminology is standard -- but it is
certainly natural for our purposes:

\begin{defi}
Le $X$ be a simplicial set. The \emph{$0$-skeletal fundamental groupoid}
of $X$ is the full subcategory of the fundamental groupoid of the
the geometric realisation $|X|$ whose object set is the
image of the $0$-skeleton.
In other words, we only allow vertices as endpoints for paths.

Let $\CG$ be a groupoid. A \emph{simplicial $K(\CG,1)$} is a
simplicial set $X$, together with a bijection between the object set of $\CG$
and the $0$-skeleton of $X$, inducing an isomorphism between
the $\CG$ and the $0$-skeletal fundamental groupoid of $X$,
and such that the higher homotopy groups of each component of $|X|$ vanish.
\end{defi}

Note that this makes sense even when $\CG$ is not connected.

\begin{prop}
The nerve $\CN\CG$ is a simplicial $K(\CG,1)$.
\end{prop}

\begin{proof}
The simplices of $\CG_{x_0} \backslash \CN\widetilde{\CG}_{x_0}$
are indexed by $\bar$ symbols
$\phantom{.}_{x}[g_1|\dots |g_n]$. It is readily seen
that the map $\CN\widetilde{\CG}_{x_0} \to \CN\CG,
g_0[g_1|\dots |g_n] \mapsto [g_1|\dots |g_n]$ induces an
isomorphism $\CG_{x_0} \backslash \CN\widetilde{\CG}_{x_0}\simeq \CN\CG$.
\end{proof}

Assume now that $\CG$ is the Garside groupoid of a Garside category
$(\CC,\phi,\Delta)$. The Garside structure provides us with
a substitute for $\CN\CG$, the \emph{Garside nerve} $\gar\CG$.

\begin{defi}
Let $(\CC,\phi,\Delta)$ be a Garside category with Garside groupoid $\CG$.
The \emph{Garside nerve} of $\CG$ is the simplicial set 
$\CN_{\Delta} \CG$ consisting of $n$-simplices of $\CN\CG$
$$\xymatrix{
x_0 \ar[r]^{f_1} & x_1 \ar[r]^{f_2} & x_2 \ar[r]^{f_3}  & x_3  \ar@{..>}[r]
& x_{n-1} \ar[r]^{f_{n}} & x_{n}}.$$
such that $f_1\dots f_n\leq \Delta_{x_0}$ (the faces and degeneracy maps
are the same as that of $\CN\CG$).
\end{defi}

The maximal length of non-degenerate simplices
is precisely the Garside dimension of $\CC$.

\begin{theo}
\label{theokaggun}
Let $(\CC,\phi,\Delta)$ be a Garside category with Garside groupoid
$\CG$. The Garside nerve $\CN_{\Delta}\CG$ is a simplicial $K(\CG,1)$
of dimension $\dim_{\Delta}\CC$.
It is finite if and only $(\CC,\phi,\Delta)$ has finite type.
\end{theo}

The geometric realisation of a simplicial set $X$ is obtained by
glueing geometric simplicial corresponding to the \emph{non-degenerate}
simplices of $X$. 

\begin{coro}
\label{corokaggun}
Let $(\CC,\phi,\Delta)$ be a Garside category with Garside groupoid
$\CG$. Let $x_0$ be an object of $\CG$.
Let $\Gamma$ be the graph whose vertex set is 
$\CG(x_0,-)$,
two vertices $f,g\in \CG(x_0,-)$ being connected by an edge
when $f\neq g$ and either $fg^{-1}$ or $gf^{-1}$ is simple.
The realisation of the simplicial complex $\Flag(\Gamma)$ is contractible.
\end{coro}

\begin{proof}
Let $\{f_0,\dots,f_n\}$ be a simplex in $\Flag(\Gamma)$.
For any $i\neq j$, we have either $f_i^{-1}f_j$ or $f_j^{-1}f_i$.
Both cannot happen simultaneously, because the product of two non-trivial
simples cannot be $1$. Write $f_i < f_j$ if $f_i^{-1}f_j$ and
$f_j < f_i$ otherwise. Let $\leq$ be the reflexive closure of $<$.

The relation $\leq$ is transitive:
let $i,j,k$ be such that $f_i< f_j$ and $f_j< f_k$;
we have either $f_i< f_k$ or $f_k< f_i$; if $f_k< f_i$, then
there are simple elements $s,t,u$ such that
$f_is=f_j$, $f_jt=f_k$, $f_ku=f_i$, thus $f_istu=f_i$ and,
by cancellativity, $stu=1$, which contradicts atomicity.

It follows that $\leq$ is a total ordering of $\{f_0,\dots,f_n\}$. 
Up to reordering, we may assume that $f_i\leq f_j \Leftrightarrow i\leq j$.

It is clear that $(f_1,\dots,f_n)$ is a non-degenerate $n$-simplex of
the Garside nerve $\CN_{\widetilde{\Delta}}\widetilde{\CG}_{x_0}$. 
Conversely, any non-degenerate
$n$-simplex of
$\CN_{\widetilde{\Delta}}\widetilde{\CG}_{x_0}$ corresponds to a
unique $n$-simplex of $\Flag(\Gamma)$.

Because the face operators are compatible with this bijection,
and because the geometric realisation of a simplicial set is precisely
obtained by glueing geometric simplices associated with non-degenerate
simplices, the realisations $|\Flag(\Gamma)|$ and $|\CN_{\widetilde{\Delta}}
\widetilde{\CG}_{x_0}|$ are homotopy equivalent. Since $\widetilde{\CG}_{x_0}$
is equivalent to the trivial category, the contractibility
of $|\CN_{\widetilde{\Delta}}
\widetilde{\CG}_{x_0}|$ follows from the theorem.
\end{proof}

The corollary is very useful for topological applications, in 
particular in the proof of the $K(\pi,1)$ conjecture for complex reflection
arrangements (\cite{kapiun}): by contrast with simplicial sets,
actual simplicial complexes may appear as nerves of open coverings.

\begin{remark}
Since we have omitted the proof of \ref{theokaggun}, it is all too convenient
to present \ref{corokaggun} as a \emph{corollary} of \ref{theokaggun}.
However, when actually checking the details of the proof, one rather goes
the other way around: one first shows that $|\Flag(\Gamma)|$ is contractible
(using, for example, Bestvina's techniques),
then that $|\CN_{\widetilde{\Delta}} \widetilde{\CG}_{x_0}|$ is contractible
(using the argument that we have presented as a ``proof'' of
Corollary \ref{corokaggun}), then that $\CN_{\Delta}{\CG}$ is 
a simplicial $K(\CG,1)$ (using Galois theory).
\end{remark}

\section{Non-positively curved aspects, after Bestvina}
\label{section8}

Let
$A$ is an Artin group of finite type, with Garside element $\Delta$.
In his beautiful article \cite{bestvina}, Bestvina made the
crucial observation that $A/\left< \Delta^2 \right>$ is
very close to being a hyperbolic group. He constructed a simplicial
complex $\CX$ together with a simplicial action
 $A/\left< \Delta^2 \right>$, and equipped $\CX$ with a ``non-symmetric''
distance with non-positive curvature features.

From this, he was able to prove that any periodic element $\gamma\in A$
is conjugate to some element with canonical length one:
$$(\exists p,q\in\BZ_{> 0}, \gamma^q=\Delta^p) \Rightarrow (\exists
s \text{ simple}, \exists k\in \BZ, \text{ $\gamma$ is conjugate to
$s\Delta^k$}).$$
This result is an immediate consequence of his Theorem 4.5 and has
many consequences: e.g., that, for given $p,q$, there are only a finite
number of conjugacy classes of $\frac{p}{q}$-periodic elements.

Charney-Meier-Whittlesey have rewritten Bestvina's proof in the language
of Garside monoids, see \cite{cmw}.
As expected, there is no obstruction to
working with categories:

\begin{theo}[after Bestvina and Charney-Meier-Whittlesey]
\label{theobestvina}
Let $(\CC,\phi,\Delta)$ be a cyclic categorical Garside structure,
with Garside groupoid $\CG$.
Let $p,q\in \BZ_{> 0}$, let $\gamma\in\CG$ be a $\frac{p}{q}$-periodic
element. There exists a simple morphism $s\in \CC$ and an integer
$k\in \BZ$ such that
$$ss^{\phi^{-k}}s^{\phi^{-2k}}\dots s^{\phi^{-(q-1)k}} = \Delta$$
and such that
$\gamma$ is conjugate (in the groupoid sense)
to $s\Delta^k$.
\end{theo}

\begin{proof} Long but easy translation exercise, left to the reader
-- sections 2,3,6 from \cite{cmw} should be rewritten in categorical
language, to obtain an analog of their Corollary 6.9.
\end{proof}

\begin{remark}
\begin{itemize}
\item[(i)] For simplicity, we only consider here \emph{cyclic}
Garside categories, because when the category is not cyclic some
definitions and arguments from \cite{bestvina} and \cite{cmw} have
no clear analogs. However, it is very likely that something can be
said about the non-cyclic case.
\item[(ii)] The version for categories of Bestvina's
complex $\CX$ should be thought of as a quotient of
the Garside nerve $\CN_{\Delta}\CG$, where all vertices are written
in normal form and the $\Delta$'s at the end are ``forgotten''.
\item[(iii)] In the next sections, we will improve the above theorem
by showing that, up to changing
the Garside structure, \emph{any periodic element is conjugate to
an element with canonical length zero!}, \ie, may be viewed
as a Garside element.
\item[(iv)] 
A key step in Bestvina's approach is a Cartan fixed point theorem
(\cite[Theorem 3.15]{bestvina}) which, adapted to our setting,
implies that the action of $\gamma$ on $\CX$
``leaves a simplex [...] invariant (and fixes its barycenter).''
To prove that any periodic element is conjugate to
an element with canonical length zero, we will make use of
an algebraic operation which replaces the category $\CC$
by a \emph{divided} Garside category (see next section).
At the level of Garside nerves, this corresponds to taking
some sort of barycentric subdivision. A consequence is that
barycentres of simplices of $\CX$ become vertices in the divided
version of $\CX$ -- one then concludes using the Cartan fixed point theorem.
Because it allows for more explicit proofs, the next sections
are written in algebraic language; however,
readers with good geometric intuition should keep in mind that everything
can thought of in terms of Garside nerves and Bestvina complexes.
When adding more vertices, one increases the number of objects in the
$0$-skeletal fundamental groupoid, without changing the fundamental groups --
in algebraic terms, this will be naturally phrased in terms of an 
\emph{equivalence} of categories between $\CG$ and its divided version.
\end{itemize}
\end{remark}

\section{Divided Garside categories}
\label{section9}

We now proceed to our main construction, which
is a general procedure, starting with a Garside
category $\CC$, to obtain a family $(\CC_m)_{m\in \BZ_{\geq 1}}$
of Garside categories. When $m=1$, one recovers $\CC$.
When $\CC$ has finite type, all $\CC_m$ have finite type (but they
usually have more objects than $\CC$).
When $\CC$ is $k$-cyclic, $\CC_m$ is $mk$-cyclic.
For $m\neq n$, the categories $\CC_m$ and $\CC_n$ are usually
not isomorphic nor equivalent, but 
their Garside groupoids $\CG_m$ and $\CG_n$ are equivalent as
categories.

\begin{defi}
Let $(\CC,\phi,\Delta)$ be a Garside triple and $m\in \BZ_{\geq 1}$.
A \emph{$m$-subdivision of $\Delta$} is a sequence
$f=(f_1,f_2,\dots,f_m)$ of composable
$\CC$-morphisms such that 
$$\prod_{i=1}^{m} f_i = \Delta.$$

We denote by $D_m(\CC,\phi,\Delta)$ (or simply $D_m(\CC)$, or simply $D_m$)
the set of $m$-subdivisions of $\Delta$.
\end{defi}

By ``$\prod_{i=1}^{m} f_i = \Delta$'', we of course mean that that the $f_i$'s
are indeed composable (the target of $f_i$ is the source of $f_{i+1}$)
and that their product is $\Delta_{x_1}$, where $x_1$ is the source of $f_1$.
This implies that the target of $f_m$ is $x_1^{\phi}$.
This also implies that each $f_i$ is simple, because
one property of Garside categories stipulates that factors of simple elements
are simple.

Note that,
contrary to our earlier notations, we 
choose to label objects starting at $1$ and not $0$, according to:

$$\xymatrix{
x_1 \ar[r]^{f_1} \ar@/_1em/[rrrrr]_{\Delta}
& x_2 \ar[r]^{f_2} & x_3 \ar[r]^{f_3}  & x_4  \ar@{..>}[r]
& x_{m} \ar[r]^{f_{m}} & x_{1}^{\phi}}.$$

\begin{convention}
In the sequel, whenever $(a_1,\dots,a_m)$ is a sequence of $\CC$-objects
or $\CC$-morphisms, we extend the notation $a_i$ to all $i\in \BZ_{\geq 1}$
by recursively setting $a_{m+i}:=a_i^{\phi}$.
(Because we have assumed that $\phi$ is invertible, we may actually 
extend our index set to $\BZ\leq 0$, although we won't use it).
\end{convention}

To illustrate this convention, we observe that,
in the above commutative diagram, we may say that 
``the target of $f_i$ is $x_{i+1}$'' without worrying about
the case $i=m$.

The object set of the $m$-th divided category $\CC_m$ will be $D_m$.
We are going to define the category $\CC_m$ by means
of its germ of simples. We have to define an oriented graph $S_m$
on the vertex set $D_m$ and endow it with a partial product structure.

Let $f=(f_1,\dots,f_m)$ and $g=(g_1,\dots,g_m)$ be two 
elements of $D_m$.
An element $s\in S_{m,f{\to} g}$ is, by definition, a 
sequence $s:=(s_1,\dots,s_m)$ of $\CC$-simples,
each $s_i$ going from the source to $f_i$ to the
source of $g_i$, forming a simple commutative diagram
in $\CC$ as follows:
$$
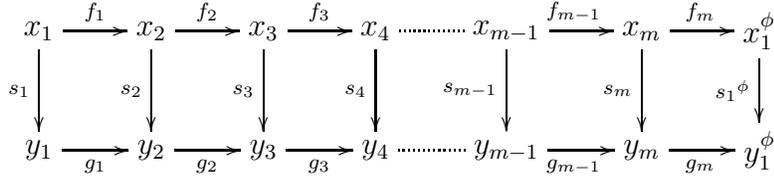
\xymatrix{
 x_1 \ar[r]^{f_1} \ar[d]_{s_1} & x_2 \ar[r]^{f_2} \ar[d]_{s_2}
& x_3 \ar[r]^{f_3} \ar[d]_{s_3} 
& x_4 \ar[d]_{s_4}
\ar@{..}[r] & x_{m-1} \ar[d]_{s_{m-1}}  
\ar[r]^{f_{m-1}} & x_m \ar[r]^{f_m} \ar[d]_{s_m} & x_1^{\phi}
\ar[d]_{{s_1}^{\phi}} \\
 y_1 \ar[r]_{g_1}  & y_2 
 \ar[r]_{g_2} & y_3  \ar[r]_{g_3} &
y_4  \ar@{..}[r] & y_{m-1}   \ar[r]_{g_{m-1}}
& y_m \ar[r]_{g_m} & y_1^{\phi} 
 }$$
and such that the above diagram may be completed by
 diagonal simple $\CC$-morphisms to obtain a commutative diagram:
$$
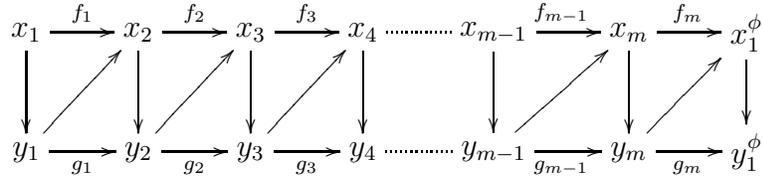
\xymatrix{
 x_1 \ar[r]^{f_1} \ar[d] & x_2 \ar[r]^{f_2} \ar[d]
& x_3 \ar[r]^{f_3} \ar[d] 
& x_4 \ar[d]
\ar@{..}[r] & x_{m-1} \ar[d]
\ar[r]^{f_{m-1}} & x_m \ar[r]^{f_m} \ar[d] & x_1^{\phi}  \ar[d] \\
 y_1 \ar[r]_{g_1} \ar[ur] & y_2  \ar[ur]
 \ar[r]_{g_2} & y_3  \ar[ur] \ar[r]_{g_3} &
y_4  \ar@{..}[r] & y_{m-1}  \ar[ur] \ar[r]_{g_{m-1}}
& y_m \ar[r]_{g_m} \ar[ur] & y_1^{\phi} 
 }$$
By cancellativy in $\CC$, there is at most one way to add the diagonal
arrows, so we may as well consider a morphism to be the whole diagram,
including those diagonal arrows.

Let us put it in other words.
A simple morphism $f\stackrel{s}{\to} g$
consists of factorisations $f_i=s_is'_i$, one 
for each $i=1,\dots,m$, such
that $g_i=s'_is_{i+1}$ for all $i\geq 1$.

Let $f,g,h\in D_m$, let $s\in S_{m,f\to g}$ 
and $t\in S_{m,g\to h}$. This corresponds to
 a commutative diagram with $\CC$-simple arrows:
$$
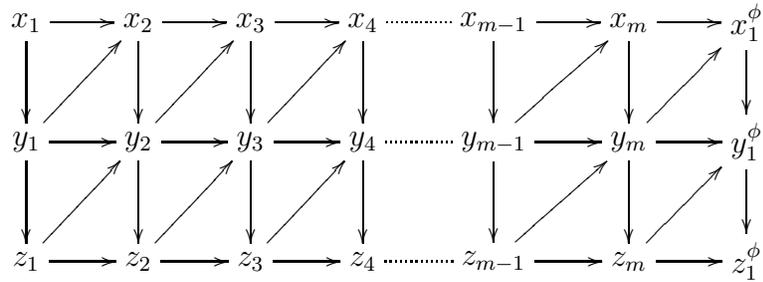
\xymatrix{
 x_1 \ar[r] \ar[d] & x_2 \ar[r] \ar[d]
& x_3 \ar[r] \ar[d] 
& x_4 \ar[d]
\ar@{..}[r] & x_{m-1} \ar[d]
\ar[r] & x_m \ar[r] \ar[d] & x_1^{\phi}  \ar[d] \\
 y_1 \ar[r] \ar[ur] \ar[d] & y_2 \ar[r] \ar[ur] \ar[d]
& y_3 \ar[r] \ar[ur] \ar[d] 
& y_4 \ar[d] 
\ar@{..}[r] & y_{m-1} \ar[d] \ar[ur]
\ar[r] & y_m \ar[r] \ar[ur] \ar[d] & y_1^{\phi}  \ar[d] \\
 z_1 \ar[r] \ar[ur] & z_2  \ar[ur]
 \ar[r] & z_3  \ar[ur] \ar[r] &
z_4  \ar@{..}[r] & z_{m-1}  \ar[ur] \ar[r]
& z_m \ar[r] \ar[ur] & z_1^{\phi} 
 }$$

We say that $s$ and $t$ are \emph{compatible} if the diagram may
be completed by simple arrows in $\CC$ to obtain a
commutative diagram:
$$
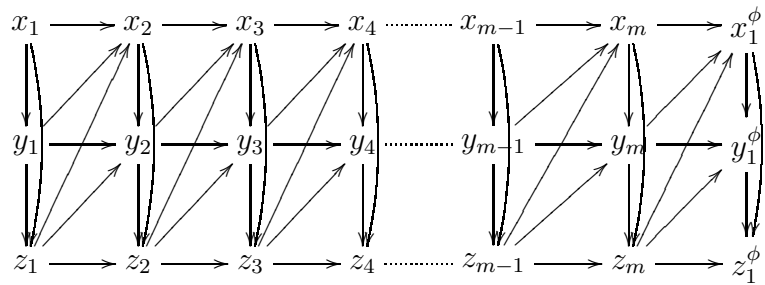
\xymatrix{
 x_1 \ar[r] \ar[d] \ar@/^0.5em/[dd] & x_2 \ar[r] \ar[d]  \ar@/^0.5em/[dd]
& x_3 \ar[r] \ar[d]  \ar@/^0.5em/[dd]
& x_4 \ar[d]  \ar@/^0.5em/[dd]
\ar@{..}[r] & x_{m-1} \ar[d]  \ar@/^0.5em/[dd]
\ar[r] & x_m \ar[r] \ar[d]  \ar@/^0.5em/[dd]
 & x_1^{\phi}  \ar[d]   \ar@/^0.5em/[dd] \\
 y_1 \ar[r] \ar[ur] \ar[d] & y_2 \ar[r] \ar[ur] \ar[d]
& y_3 \ar[r] \ar[ur] \ar[d] 
& y_4 \ar[d] 
\ar@{..}[r] & y_{m-1} \ar[d] \ar[ur]
\ar[r] & y_m \ar[r] \ar[ur] \ar[d] & y_1^{\phi}  \ar[d] \\
 z_1 \ar[r] \ar[ur] \ar[uur] & z_2  \ar[ur] \ar[uur]
 \ar[r] & z_3  \ar[ur] \ar[uur] \ar[r] &
z_4  \ar@{..}[r] & z_{m-1}  \ar[ur] \ar[uur] \ar[r]
& z_m \ar[r] \ar[ur] \ar[uur] & z_1^{\phi} 
 }$$

In other words, $s=(s_1,\dots,s_m)$ and $(t_1,\dots,t_m)$
are compatible if and only if each $s_it_i$ may be multiplied in $S$ and
$$(s_1t_1,\dots,s_mt_m) \in S_m.$$
We take as partial product structure on $S_m$ the
map sending compatible $s$ and $t$ as above to
$$st:=(s_1t_1,\dots,s_mt_m).$$
This defines a germ structure $$\CS_m$$
with underlying graph $S_m$.

Given an object 
$$f=(\xymatrix{
x_1 \ar[r]^{f_1} & x_2 \ar[r]^{f_2} & x_3 \ar[r]^{f_3}  &   \ar@{..}[r]
& x_{m-1} \ar[r]^{f_{m-1}} & x_m \ar[r]^{f_m} & x_1^{\phi}})$$
in $D_m(\CC)$, the poset
$(S_m(\CC)_{f\to},\leq)$ is isomorphic to
$$\prod_{i=1}^m ([1_{x_i},f_i],\leq),$$
where $[1_{x_i},f_i]$ denotes the interval between $1_{x_i}$ and
$f_i$ in $(S_{x_i\to},\leq)$.
Since each $(S_{x_i},\leq)$ is a lattice, $(S_m(\CC)_{f\to},\leq)$
is a lattice, whose maximal element is
$$\xymatrix{
 x_1 \ar[r]^{f_1} \ar[d]_{f_1} & x_2 \ar[r]^{f_2} \ar[d]_{f_2}
& x_3 \ar[r]^{f_3} \ar[d]_{f_3} 
& x_4 \ar[d]_{f_4}
\ar@{..}[r] & x_{m-1} \ar[d]_{f_{m-1}}
\ar[r]^{f_{m-1}} & x_m \ar[r]^{f_m} \ar[d]_{f_m} & x_1^{\phi}
 \ar[d]_{f_1^{\phi}} \\
 x_2 \ar[r]_{f_2} & x_3  
 \ar[r]_{f_3} & x_4 \ar[r]_{f_4} &
x_5  \ar@{..}[r] & x_m \ar[r]_{f_m}
& x_1^{\phi} \ar[r]_{f_1^{\phi}}  & x_2^{\phi} 
 }$$
(which is to be completed by diagonal trivial $\CC$-morphisms $1_{x_i}$).

We denote by $\Delta_{m,f}$ this element of $S_m$. We define
an automorphism $\phi_m$ of $S_m$ acting on vertices by
$$(f_1,f_2,\dots,f_m)^{\phi_m}:=(f_2,\dots,f_m,f_1^{\phi})$$
and on arrows by
$$(s_1,s_2,\dots,s_m)^{\phi_m}:=(s_2,\dots,s_m,s_1^{\phi}).$$
This is compatible with the partial product and induces an automorphism
$\phi_m$ of $\CS_m$, such that $\Delta_m$ is a natural transformation
from the identify functor to $\phi_m$.

\begin{defi}
The \emph{$m$-divided category} associated with $(\CC,\phi,\Delta)$
is the free category $\CC_m$ on the germ $\CS_m$.
The groupoid of fractions of $\CC_m$ is denoted by $\CG_m$.
\end{defi}

\begin{theo}
The triple 
$(\CC_m, \phi_m,\Delta_m)$ is a categorical Garside structure.
\end{theo}

\begin{proof}
It is clear from the preceding discussion.
\end{proof}

Let $(\CC,\phi,\Delta)$ be a Garside triple, let $m\geq 1$.

For any object $x$ of $\CC$, we consider the object of $\CC_m$ 
defined by
$$\Theta_m(x) :=  \left(
\xymatrix@1{
x \ar[r]^{1_x} 
& x \ar[r]^{1_x} & x  \ar@{..>}[r]   & x \ar[r]^{1_x}
& x \ar[r]^{\Delta} & x^{\phi}} \right).$$
If $x\stackrel{s}{\to}y$ is a simple $\CC$-morphism, we
consider the $\CC_m$-morphism 
$$\Theta_m(x) \stackrel{\Theta_m(s)}{\to} \Theta_m(y)$$
defined as the composition from top to bottom of the following
simple morphisms (note that $\Theta_m(x)$ itself is not simple):
$$
\xymatrix{
x \ar[r]^{1} \ar[d]_1
& x \ar[r]^{1} \ar[d]_1 & x \ar[d]_1 
\ar@{..>}[r]   & x \ar[r]^{1} \ar[d]_1
& x \ar[r]^{\Delta} \ar[d]_s & x^{\phi}  \ar[d]_1 \\
x \ar[r]^{1}  \ar[d]_1
& x \ar[r]^{1} \ar[d]_1 & x  \ar[d]_1 \ar@{..>}[r] 
 & x \ar[r]^{s} \ar[d]_s
& y \ar[r]^{\overline{s}} \ar[d]_1 & x^{\phi} \ar[d]_1\\
x \ar@{..>}[d] \ar[r]^{1} 
& x \ar@{..>}[d] \ar[r]^{1} & x \ar@{..>}[d] \ar@{..>}[r]   & y \ar[r]^{1}
\ar@{..>}[d]
& y \ar[r]^{\overline{s}}  \ar@{..>}[d] & x^{\phi} \ar@{..>}[d] \\
x \ar[r]^{1}  \ar[d]_1
& x \ar[r]^{s} \ar[d]_s & y  \ar[d]_1 \ar@{..>}[r] 
 & y \ar[r]^{1} \ar[d]_1
& y \ar[r]^{\overline{s}} \ar[d]_1 & x^{\phi} \ar[d]_1 \\
x \ar[r]^{s}  \ar[d]_s
& y \ar[r]^{1} \ar[d]_1 & y  \ar[d]_1 \ar@{..>}[r] 
 & y \ar[r]^{1} \ar[d]_1
& y \ar[r]^{\overline{s}} \ar[d]_1 & x^{\phi} \ar[d]_{s^{\phi}} \\
y \ar[r]^{1}  
& y \ar[r]^{1}  & y   \ar@{..>}[r] 
 & y \ar[r]^{1} & y \ar[r]^{\Delta}  & y^{\phi} 
 }
$$

\begin{theo}
\label{theo94}
The map $x\mapsto \Theta_m(x),
\left( x \stackrel{s}{\to} y \right) \mapsto 
\left( \Theta_m(x) \stackrel{\Theta_m(s)}{\longrightarrow} \Theta_m(y) \right)$
extends to a unique functor $$\Theta_m:\CC \to \CC_m$$
whose induced functor $$\CG \to \CG_m$$ is an equivalence of
categories.
\end{theo}

\begin{proof}
To check that the functor $\Theta_m:\CC \to \CC_m$ is well-defined,
one has to check that, whenever $st=u$ holds
in the germ of simples of $\CC$, one has $\Theta_m(s)\Theta_m(t)=\Theta_m(u)$
in $\CC_m$.
This is a straighforward computation.

By theorem \ref{theokaggun}, we know that a Garside group is the
fundamental groupoid of its Garside nerve (with respect to the $0$-skeleton).

A $(k-1)$-simplex of the Garside nerve of $\CG_m$ consists of the following
data:
\begin{itemize}
\item a $\CC_m$-object
$$f = \left(\xymatrix@1{
x_1 \ar[r]^{f_1} \ar@/_1em/[rrrrr]_{\Delta_{x_1}}
& x_2 \ar[r]^{f_2} & x_3 \ar[r]^{f_3}  & x_4  \ar@{..>}[r]
& x_{m} \ar[r]^{f_{m}} & x_{1}^{\phi}} \right),$$
which we call the \emph{basepoint} of the simplex (and we say
that the basepoint \emph{starts at $x_1$}),
\item a totally ordered (for $\leq$) sequence of $k$
simple $\CC_m$-morphisms
with source $f$ or, equivalently,
 a totally ordered (for $\leq$) sequence of $k+1$
simple $\CC_m$-morphisms
with source $f$ and whose last term is $\Delta_{m}$ or, again equivalently,
a factorisation 
$f_i = f_{i,1}\dots f_{i,k+1}$ of each $f_i$ into $k+1$ simple $\CC$-morphisms
-- in the latter description, the successive
simple $\CC_m$-morphisms from $f$ associated to the factorisations
are
$$\xymatrix{
 \cdot \ar[rr]^{f_1} \ar[d]_{u_{1,j}} & & \cdot
\ar[rr]^{f_2} \ar[d]_{u_{2,j}}
& & \cdot \ar[rr]^{f_3} \ar[d]_{u_{3,j}}
& & \cdot \ar[d]_{u_{4,j}}
\ar@{..}[rr] & & \cdot \ar[d]_{u_{m-1,j}}
\ar[rr]^{f_{m-1}} & & \cdot \ar[rr]^{f_m} \ar[d]_{u_{m,j}}
 & & \cdot  \ar[d]_{u_{1,j}^{\phi}} \\
 \cdot \ar[rr]_{v_{1,j}u_{2,j}} & & \cdot \ar[rr]_{v_{2,j}u_{3,j}}
& & \cdot \ar[rr]_{v_{3,j}u_{4,j}}
& & \cdot \ar@{..}[rr] & & \cdot
\ar[rr]_{v_{m-1,j}u_{m,j}} &  & \cdot
\ar[rr]_{ v_{m,j}u_{1,j}^{\phi} } & & \cdot },$$
where $u_{i,j}:=f_{i,1}\dots f_{i,j}$ and
$v_{i,j}:=f_{i,j+1}\dots f_{i,k+1}$.
\end{itemize}

Since $f_1\dots f_m= \Delta_{x_1}$, we see that a $(k-1)$-simplex
of $\CN_{\Delta_m}\CG_m$ with basepoint starting at $x_1$ is nothing
but a factorisation of $\Delta_{x_1}$ into $km$ terms, that is,
a $(mk-1)$-simplex of $\CN_{\Delta}\CG$ with basepoint $x_1$. 

In other words, $\CN_{\Delta_m}\CG_m$ is the $m$-th edgewise subdivision
of $\CN_{\Delta}\CG$, in the sense of B\" okstedt-Hsiang-Madsen (see
section 1 of \cite{bhm}). By \cite[Lemma 1.1]{bhm}, the realisations
$|\CN_{\Delta_m}\CG_m|$ and $|\CN_{\Delta}\CG|$ are homeomorphic.
This implies that there exists an equivalence of categories between
their $0$-skeletal fundamental groupoids, which are precisely
$\CG_m$ and $\CG$.

The last thing to check is that $\Theta_m$ actually provides such an
equivalence of categories. B\" okstedt-Hsiang-Madsen's definition of an
homeomorphism\footnote{B\" okstedt-Hsiang-Madsen's homeomorphism $D_m$,
used in this paragraph, should not be confused with our $D_m$.} $D_m$
between $|\CN_{\Delta_m}\CG_m|$ and $|\CN_{\Delta}\CG|$
is constructive and, to conclude,
it suffices to check that $D_m(\Theta_m(x))=x$ for all $\CC$-object
(and that $D_m$ behaves as expected with respect to the $1$-skeleton).
A $0$-simplex $(f_1,\dots,f_m)$ in $\CN_{\Delta_m}\CG_m$ with basepoint $x_1$
corresponds to a
$(m-1)$-simplex in $\CN_{\Delta}\CG$ whose vertices are the sources
of the $f_i$'s and the edges are partial products of the first $m-1$ $f_i$'s.
By definition, $D_m(f_1,\dots,f_m)$ is taken to be the barycentre of
that $(m-1)$-simplex. But the $(m-1)$-simplex associated with $\Theta_m(x)$
is completely degenerate,
all its vertices are $x_1$ and all its edges are $1_{x_1}$.
Therefore 
it collapses to the single point $x_1$ in the realisation. 
\end{proof}

Let us conclude this section with some easy basic properties.

\begin{prop}
Let $(\CC,\phi,\Delta)$ be a finite type 
categorical Garside structure of Garside dimension $n$.
There exists a polynomial $Z{(\CC,\phi,\Delta)}$, of degree
at most $n$ and with integral coefficients, such that for all $m\geq 1$
the number of elements in $D_m$ is $Z{(\CC,\phi,\Delta)}(m)$
\end{prop}

\begin{proof}
This follows from the fact that the number of $m$-chains
in a finite poset is a polynomial in $m$.
\end{proof}

\begin{prop}
Let $(\CC,\phi,\Delta)$
be a categorical Garside structure of Garside dimension $n$.
Then $$\dim_{\Delta_m} \CC_m = \dim_{\Delta} \CC \qquad \text{and} \qquad
\dim_{\Delta_m} \CC_m^{\phi_m} \leq \frac{\dim_{\Delta} \CC^{\phi}}{m}.$$
\end{prop}

\begin{prop}
Let $(\CC,\phi,\Delta)$ be a categorical Garside structure.
Let $p,q,e$ be positive integers. The natural bijection
$$D_{eq}(\CC,\phi,\Delta) \simeq D_{e} (\CC_q,\phi_q,\Delta_q)$$ 
induces an isomorphism
$$ \CC_{eq}^{\phi_{eq}^{ep}} \simeq (\CC_q^{\phi_q^p})_e.$$
\end{prop}

\section{Periodic elements are Garside elements}
\label{section10}

Recall that a $\frac{p}{q}$-periodic element in a Garside groupoid is a loop
$\gamma$ such that $$\gamma^q= \Delta^p.$$
After Bestvina, we observed in Theorem \ref{theobestvina} that,
if the Garside category is cyclic, any periodic element is
conjugate to some $s\Delta^k$, where $s$ is a simple element such
that $$ss^{\phi^{-k}}s^{\phi^{-2k}}\dots s^{\phi^{-(q-1)k}} = \Delta.$$
To alleviate notations, we set
$s_i:=s^{\phi^{-(i-1)k}}$ and $\underline{s}:=(s_1,\dots, s_q)$. The above
identity expresses that $\underline{s}$ is an object of $\CC_q$.
Using $1\stackrel{\phi}{\Rightarrow}\Delta$, we obtain
$$\Delta^p = (s\Delta)^q = s_1\dots s_q\Delta^{qk} =\Delta^{qk+1}$$
thus $$p=qk+1.$$
\begin{theo}
\label{theo101}
Let $(\CC,\Delta,\phi)$ be a cyclic categorical Garside structure
with associated groupoid $\CG$. As above,
let $\rho=s\Delta^k$ be a $\frac{p}{q}$-periodic
loop in $\CG$. Consider the $q$-divided Garside groupoid $\CG_q$
and the functor $\Theta_q:\CG \to \CG_q$.
Then $\Theta_q(\rho)$ is conjugate (in the groupoid sense) to
the element $\Delta_q^{p}$ (product of $p$ successive
Garside elements of $\CG_q$) with source
$\underline{s} := (s_1,\dots,s_q)$.
\end{theo}

This theorem should be thought of as an algebraic
Ker\'ek\'art\'o-Brouwer-Eilenberg
theorem for Garside categories -- this will become clearer in the next
section, as we rephrase this in terms of $S^1$-spaces.
It gives a positive answer to the 
categorical rephrasing of Question 4.

\begin{proof}
Consider the configuration space $U_{q}$ of
$q$ unordered distinct points on the unit circle $S^1$
(or, more intuitively, ``beads on a necklace'').
Consider the subset $U_{q,m}$ consisting of configurations
in $S^1-\mu_m$, where $\mu_m$ is the group of $m$-th roots
of unity.
For $i=1,\dots,m$, let $V_i$ be the connected component of $S^1-\mu_m$
consisting of points with argument in
 $(\frac{2\pi(i-1)}{m},\frac{2\pi i}{m})$.

 Because the connected components of $U_{q,m}$ are contractible,
we may consider the fundamental groupoid  $NB_{q,m}$ of
$U_{q}$ with respect to these components:
\begin{itemize}
\item objects of $NB_{q,m}$ are connected components
of $U_{q,m}$; each component is uniquely determined by the
sequence $(n_1,\dots,n_m)$ where $n_i$ is the number of beads in $V_i$;
conversely, any sequence $(n_1,\dots,n_m)\in(\BZ_{\geq 0})^m$ with sum $q$
corresponds to an object;
\item if $C$ and $C'$ are objects,
$\Hom_{NB_{q,m}}(C,C') := \varinjlim_{x\in C}
\varinjlim_{x'\in C'} \Hom_{\pi_1(U_{q})}(x,x').$
\end{itemize}
Elements of $NB_{q,m}$ are \emph{necklace braids with $q$ beads and
$m$ sectors}, or simply \emph{necklace braids}. From now on, 
we identify objects of $NB_{q,m}$ with their associated sequences.

Let $i\in\{1,\dots,m\}$ and consider a component $(n_1,\dots,n_m)$.
If $n_i>0$, we define the \emph{$i$-slide with source $(n_1,\dots,n_m)$}
as the element of $NB_{q,m}$ obtained by ``sliding'' a single 
bead from $V_i$ to $V_{i-1}$ ($V_m$ if
$i=1$) by decreasing its argument and crossing once the point with argument
$\frac{2\pi(i-1)}{m}$.
We denote this element by $\sigma_i$, regardless of its source (in
a given formula, the symbol $\sigma_i$ should be interpreted as the
\emph{only} possible $i$-slide whose source is as provided by the context).

It is clear that the necklace
braid groupoid is generated by all slides and inverses of slides.
One may easily write a presentation (whenever it makes sense, slides
commute).

Let $A$ be an alphabet, together with a permutation
$A\to A, a\mapsto a^{\phi}$. Let $W_{q,m}$ the set of $m$-tuples
$w=(w_1,\dots,w_m)$ of words in $A^*$
whose concatenation $w_1\dots w_m$ has length $q$.

Let $w=(w_1,\dots,w_m)\in W_{q,m}$
We say that a necklace braid $\beta$ is \emph{compatible} with $w$
if its source is $(l(w_1),\dots,l(w_m))$.
Assume that $\sigma_i$ is compatible with $w$, \ie, that
we may write $w_i=aw'_i$ (with $a\in A$
and $w'_i\in A^*$). We define an element $w\cdot \sigma_i\in W_{q,m}$
as follows:
\begin{itemize}
\item if $i>1$, we set $w\cdot \sigma_i := (w_1,\dots,w_{i-2},w_{i-1}a,w'_i,
w_{i+1},\dots, w_m)$,
\item if $i=1$, we set 
$w\cdot \sigma_i := (w'_1,w_2,\dots, w_{m-1},w_ma^{\phi})$.
\end{itemize}

This extends to a right action of $NB_{q,m}$ on
$W_{q,m}$. By this, we mean that there is a category
$W_{q,m} \circlearrowright NB_{q,m}$ whose object set is
$W_{q,m}$, such that $\Hom_{W_{q,m} \circlearrowright NB_{q,m}}(w,-)$
is the set of necklace braids compatible with $w$, and such that,
if as above $w$ and $\sigma_i$ are compatible, the corresponding morphism
with source $w$ has target $w\cdot \sigma_i$.

Pursuing our trend of convenient abusive notation,
when the source in $W_{q,m}$ is specified,
we denote by $\sigma_i$ the only
possible $W_{q,m} \circlearrowright NB_{q,m}$-morphism with this source
and that is associated
with an $i$-slide (there is \emph{at most} one such morphism).

To prove the theorem, we apply this to $q=m$ and $A:=\{s^{\phi^k} | k\in \BZ\}$
and we only look at $O$, the $NB_{q,q}$ orbit of $(\varepsilon,\dots,
\varepsilon, s_1 \dots s_q)$, and the corresponding subcategory
$O \circlearrowright NB_{q,q}$.
By evaluating each element of $A^*$ to its product in $\CC$, we obtain
a map $\psi$ from $O$ to the object set of $\CC_q$. If $\sigma_i$ is a slide
and $(w_1,\dots,w_q) \stackrel{\sigma_i}{\to} (w'_1,\dots,w'_q)$ is
a morphism in $O \circlearrowright NB_{q,q}$,
we define $\psi((w_1,\dots,w_q) \stackrel{\sigma_i}{\to} (w'_1,\dots,w'_q))$
to be the $\CC_q$-simple morphism from $(f_1,\dots,f_q):=\psi((w_1,\dots,w_q))$
to  $(f'_1,\dots,f'_q):=\psi((w'_1,\dots,w'_q))$ corresponding to the
diagram:
$$\xymatrix{
 \cdot  \ar[r]^{f_1} \ar[d]_{1} & \cdot \ar@{..}[r] \ar[d]_{1}
& \cdot \ar[r]^{f_{i-1}} \ar[d]_{1} 
& \cdot \ar[d]_{a}
\ar[r]^{f_i} & \cdot \ar[d]_{1}  
\ar@{..}[r] & \cdot \ar[r]^{f_q} \ar[d]_{1} & \cdot
\ar[d]_{1} \\
 \cdot \ar[r]_{f'_1}  & \cdot 
 \ar@{..}[r] &\cdot  \ar[r]_{f'_{i-1}} &
\cdot  \ar[r]_{f'_i} & \cdot  \ar@{..}[r]
& \cdot \ar[r]_{f'_q} & \cdot 
 }$$
where $a$ is the first letter of $w_i$.
This extends to a functor $$\psi:O \circlearrowright NB_{q,q}
\to \CC_q.$$

If $w\in O$ is of the form $(\varepsilon,\dots,\varepsilon,w_n)$,
then  $\Delta_{q,\psi(w)} = \psi((\sigma_{n}\sigma_{n-1}\dots \sigma_1)^q)$.
If $w\in O$ is of the form $(a_1,\dots,a_q)$, where $a_1,\dots,a_q$ are
letters, then
 $\Delta_{q,\psi(w)} = \psi(\sigma_{n}\sigma_{n-1}\dots \sigma_1)$.

Consider the following elements of $O$:
$$w:=(\varepsilon,\dots,\varepsilon,s_1\dots s_q).$$
$$w':=(s_1,\dots,s_q)$$
and the following words in the alphabet $\{\sigma_1,\dots,\sigma_q\}$:
$$\beta_1:= \sigma_{n}\sigma_{n-1}\dots \sigma_1 $$
$$\beta_2:= (\sigma_{n}\sigma_{n-1}\dots \sigma_2)(\sigma_n\sigma_{n-1}\dots
\sigma_3)(\sigma_{n}\sigma_{n-1}\dots \sigma_4) \dots (
\sigma_{n}\sigma_{n-1})(\sigma_{n}).$$

In the category $NB_{q,q}$, one easily checks (make a picture!)
the relation
$$\beta_1^{qk+1} \beta_2 = \beta_2 \beta_1^{qk+1},$$
where both sides are first expanded as words in the $\sigma_i$'s
and then interpreted as morphisms with source $(0,\dots,0,q)$.
As a consequence, we obtain in $O\circlearrowright NB_{q,q}$
the relation
$$\beta_1^{qk+1} \beta_2 = \beta_2 \beta_1^{qk+1},$$
where both sides are now interpreted as morphisms with source
$w$. By functoriality, we obtain in $\CC_q$ the relation
$$\psi(\beta_1^{qk+1} \beta_2) = \psi(\beta_2 \beta_1^{qk+1}).$$

From the definition of $\Theta_q$, it is clear that
$$\Theta_q(s\Delta^k) = \psi(\beta_1^{qk+1})$$
($\beta_1^{qk+1}$ interpreted here with source $w$).

We note that $\psi(w')$ coincides with the $\CC_q$-object $\underline{s}$
from the
theorem's statement. When interpreting $\beta_1$ with source
$w'$, it maps by $\psi$ to the Garside element $\Delta_{q,\underline{s}}$
with source $\underline{s}$. When interpretating $\beta_2$ with
source $w$, its target is $w'$ and it maps by $\psi$ to a certain
$\CC_q$-morphism from $\Theta_q(x)$ (where $x$ is the source of $s$)
to $\underline{s}$. The relation
$\psi(\beta_1^{qk+1} \beta_2) = \psi(\beta_2 \beta_1^{qk+1})$ may
be interpreted as expressing a conjugacy relation relation
between $\Theta_q(s\Delta^k)$ and the
$\Delta_q^p=\Delta_q^{qk+1}$ with source $\underline{s}$,
proving the theorem.
\end{proof}

One should not be surprised by the fact that, in the theorem,
 $\Theta_q(\rho)$ is conjugate to some $\Delta_q^p$ that is by definition
a $p$-periodic
(and not $\frac{p}{q}$-periodic) element in $\CG_q$. The explanation
is that $\Delta_q$, the Garside element of $\CG_q$, behaves very much
like a ``$q$-th root'' of $\Delta$.

Another potentially disturbing fact is that not all elements of 
$\CG_q$ of the form $\Delta_q^p$ are periodic: in our definition of
periodic elements, we have required them to be \emph{loops}.
The  $\Delta_q^p$ with source $\underline{s}$ appearing in
the theorem is a loop. There is actually a very simple criterion
to decide whether such a given $\Delta_q^p$ is a loop:

\begin{lemma}
\label{lemmaphipq}
Let $f=(f_1,\dots,f_q)$ be an object of $\CG_q$.
The following conditions are equivalent:
\begin{itemize}
\item[(i)] The $\CG_q$-morphism
$\Delta_q^p$ with source $f$ is a loop, thus a $p$-periodic element.
\item[(ii)] Its source
$f$ is an object of $\CG_q^{\phi_q^p}$ (the invariant subcategory
of $\CG_q$
for the $p$-th power of the diagram automorphism $\phi_q$).
\item[(iii)] One has, for all $i=2,\dots,q$, $f_i=f_1^{\phi^{-(i-1)k}}$.
\end{itemize}
\end{lemma}

\begin{proof}
Because $1\stackrel{\Delta_q}{\Rightarrow} \phi_q$,
the target of the $\Delta_q^p$ with
source $f$ is $f^{\phi_q^p}$. This shows the equivalence between (i)
and (ii).

To check the equivalence with (iii), one uses the relation
 $p=qk+1$: a direct computation shows that the target of
the $\Delta_q^p$ with source $f$ is $(f_2^{\phi^k},f_3^{\phi^{k}},\dots,
f_q^{\phi^k},f_1^{\phi^{k+1}}$. The element is a loop
if and only if, for all $i=2,\dots,q$, $f_i=f_1^{\phi^{-(i-1)k}}$.
(Note that an additional relation seem to be required, namely that
$f_1^{\phi^{k+1}}=f_q$, but as a consequence of the other relations it
rewrites as $f_1^{\phi^p}=f_1$ and
comes for free using
$1\stackrel{\Delta}{\Rightarrow} \phi$.)
\end{proof}

As suggested by the above lemma, Theorem \ref{theo101} should
be understood as part of a deeper dictionnary
between the conjugacy category of $\frac{p}{q}$-periodic elements
in $\CG$ and the fixed subcategory $\CG_q^{\phi_q^p}$.
This has many consequences. For example, the categorical
rephrasing of Question 5 admits a positive answer:

\begin{coro}
The centraliser of a periodic element in a cyclic Garside groupoid is 
a weak Garside group.
\end{coro}

\begin{proof}
We apply Theorem \ref{theo101}. Because of the equivalence of categories,
it suffices to show that the centraliser of a periodic
 power of a Garside element in a cyclic Garside groupoid is
a weak Garside group. This has been done in Corollary \ref{corocommute}.
\end{proof}

The theorem also yields a precise criterion to test Question 2:
\begin{coro}
Let $(\CC,\Delta,\phi)$ be a cyclic categorical Garside structure
with associated groupoid $\CG$. Fix positive integers $p,q$.

Let $s$ be a $\CC$-simple such that $s\Delta^k$ is a
$\frac{p}{q}$-periodic loop (as we have seen, this forces
$p=qk+1$).
The $\CG_q$-object $\underline{s}$ (as in Theorem \ref{theo101}) is
an object of $\CG_q^{\phi_q^p}$.

If $s$ and $s'$ are  $\CC$-simples such that
$s\Delta^k$ and $s'\Delta^k$ are conjugate $\frac{p}{q}$-periodic
loops, then $\underline{s}$ and $\underline{s'}$ lie in the
same connected component of $\CG_q^{\phi_q^p}$.

In particular, we have a well-defined map from the set of conjugacy
classes of $\frac{p}{q}$-periodic loops in $\CG$ to the set 
of connected components of $\CG_q^{\phi_q^p}$ (sending $\rho$
to the connected component of $\underline{s}$, where $s$ is chosen
such that $s\Delta^k$ is a summit in the conjugacy class of $\rho$).

This map is a bijection.
\end{coro}

\begin{remark}
When $p=q=1$, the corollary describes conjugacy classes of Garside elements.
\end{remark}

\begin{proof}
For any $s$ such that $s\Delta^k$ is $\frac{p}{q}$-periodic, denote
by $\nabla(s)$ the $\Delta_q^p$ with source $\underline{s}$ (see
Theorem \ref{theo101}).

Because $\nabla(s)$ is a loop, $\underline{s}$ is $\phi_q^p$-invariant
(this is Lemma \ref{lemmaphipq}).

Suppose that $s\Delta^k$ and $s'\Delta^k$ are conjugate in $\CG$. Using
Theorem \ref{theo101} and the functoriality
of $\Theta$, one sees that $\nabla(s)$ and $\nabla(s')$ are
 conjugate in $\CG_q$. Let
$c\in \Hom_{\CG_q}(\underline{s},\underline{s'})$
be such that $$\nabla(s) c =c \nabla(s').$$
Because $\nabla(s)$ and $\nabla(s')$ are products of $p$
successive $\Delta_q$'s, we deduce from
$1\stackrel{\Delta_q}{\Rightarrow} \phi_q$ that
$$\nabla(s) c^{\phi_q^p} =c \nabla(s').$$
By cancellativity, $c=c^{\phi_q^p}$. In particular, $c$ connects
$\underline{s}$ and $\underline{s}'$ in the category $\CG_q^{\phi_q^p}$.

Conversely, suppose that $\underline{s}$ and $\underline{s}'$ lie
in the same component of $\CG_q^{\phi_q^p}$.
Let $$c\in
\Hom_{\CG_q^{\phi_q^p}}(\underline{s},\underline{s}').$$
We may view $c$ as a $q$-tuple $(c_1,\dots,c_q)$
of $\CG$-morphisms such that the
diagram
$$\xymatrix{
 \cdot \ar[r]^{s_1} \ar[d]_{c_1} & \cdot \ar[r]^{s_2} \ar[d]_{c_2}
& \cdot \ar[r]^{s_3} \ar[d]_{c_3} 
& \cdot \ar[d]
\ar@{..}[r] & \cdot \ar[d]  
\ar[r]^{s_{q-1}} & \cdot \ar[r]^{s_q} \ar[d]_{c_q} & \cdot
\ar[d]_{c_1^{\phi}} \\
 \cdot \ar[r]_{s'_1}  & \cdot 
 \ar[r]_{s'_2} & \cdot  \ar[r]_{s'_3} &
\cdot  \ar@{..}[r] & \cdot  \ar[r]_{s'_{q-1}}
& \cdot \ar[r]_{s'_q} & \cdot 
 }$$
is commutative. Because $c^{\phi_q^p}=c$, we have $c_2=c_1^{\phi^{-k}}$.
In $\CG$, we have the relation
 $$s \Delta^k c_1 = s_1\Delta^k c_1= s_1 c_2 \Delta^k = c_1 s'_1\Delta^k
=c_1 s'\Delta^k.$$
We have proved that $s\Delta^k$ and $s'\Delta^k$ are conjugate.
\end{proof}

\section{Cyclic structure on the Garside nerve}
\label{section11}

\emph{This section is a sketch. It explains the connection between
the algebraic version (Theorem \ref{theo101}) and the geometric version
(Theorem \ref{mykerekjarto}) of the ``Ker\'ekj\'art\'o principle'' for
Garside categories. It
may be skipped in a first reading.}

We study a special feature of the simplicial
structure of the Garside nerve $\CN_{\Delta}\CG$ of a Garside groupoid.
When the Garside structure is $1$-cyclic, we show that the Garside nerve
is a \emph{cyclic set}, in the sense of Connes. When the Garside structure is
only $k$-cyclic, then the Garside
nerve is very close to being a cyclic set and its realisation
may be equipped with a natural structure of $S^1$-space.

In addition to the $n+1$-degeneracy maps $$s_0,\dots,s_n:
(\CN_{\Delta}\CG)_n \to (\CN_{\Delta}\CG)_{n+1}$$
consisting of inserting identity morphisms at any of the $n+1$
objects of the sequence 
$$\xymatrix{
x_0 \ar[r]^{f_1} & x_1 \ar[r]^{f_2} & x_2 \ar[r]^{f_3}  & x_3  \ar@{..>}[r]
& x_{n-1} \ar[r]^{f_{n}} & x_{n}},$$
there is another natural way to obtain a $n+1$-simplex from a $n$-simplex:
there is a unique way of completing $(f_1,\dots,f_n)$ to a sequence
$$s_{n+1}(f_1,\dots,f_n):=(f_1,\dots,f_{n+1})$$
such that $f_1\dots f_{n+1} = \Delta$.

\begin{remark}
It is clear that the image of $s_{n+1}$ is precisely $D_{n+1}$.
\end{remark}

\begin{defi}
\label{defispecial}
We call \emph{special degeneracy operator} the degree $1$ map
$$\mathfrak{s}:\CN_{\Delta}\CG\to\CN_{\Delta}\CG$$ whose restriction
to $(\CN_{\Delta}\CG)_n$ is $s_{n+1}$.

We call \emph{first face operator} the degree $-1$ map
\begin{eqnarray*}
d_0:\CN_{\Delta}\CG & \to & \CN_{\Delta}\CG \\
(f_1,\dots,f_{n}) & \mapsto & (f_2,\dots,f_{n}). 
\end{eqnarray*}
\end{defi}

\begin{lemma}
\label{lemmaconnes}
\begin{itemize}
\item[(i)] For all $(f_1,\dots,f_n)\in D_n$,
$\mathfrak{s}d_0 (f_1,\dots,f_n) = (f_2,\dots,f_n,f_1^{\phi})$.
\item[(ii)] For all $(f_1,\dots,f_n)\in (\CN_{\Delta}\CG)_{n}$,
$(d_0\mathfrak{s})^{n+1}(f_1,\dots,f_n)=(f_1^{\phi},\dots,f_n^{\phi}).$
\end{itemize}
\end{lemma}

\begin{proof} (i) follows from $1\stackrel{\Delta}{\Rightarrow} \phi$.
(ii) is an easy consequence of (i).
\end{proof}

When $\phi$ is the identity, \ie, when $\CC$ is $1$-cyclic, statement
(ii) of the lemma
says that $\CN_{\Delta}\CG$ is a \emph{cyclic set},
in the sense of Connes, \cite{connes}.
More generally, when $\CC$ is $k$-cyclic, $\CN_{\Delta}\CG$ is a
$\Lambda_k^{\op}$-object in the category of sets, in the sense
of B\"okstedt-Hsiang-Madsen, \cite[Definition 1.5]{bhm}.

\begin{remark}
Lemma \ref{lemmaconnes} can be understood in the general setting when
$\CG$ is not cyclic, in terms of an ``helicoidal
category'' generalising Connes' cyclic category: the maps $(d_0s_{n+1})^{n+1}$
equip each $(\CN_{\Delta}\CG)_n$ with a $\BZ$-action, and
the faces and degeneracy maps are $\BZ$-equivariant.
\end{remark}

\begin{theo}
Let $(\CC,\phi,\Delta)$ be a $k$-cyclic categorical Garside structure.
Let $X:=|\CN_{\Delta} \CG|$ be the realisation of the Garside nerve
of the associated groupoid. There is a canonical
structure of $S^1$-space on $X$, with respect to which
Question 1 has a positive
answer: any periodic loop of $\pi_1(X)$ is conjugate to a rotation.
\end{theo}

\section{Example: $3$-divided category of the Artin-Tits monoid of type $A_2$}

The classical Artin-Tits monoid of type $A_2$ is a Garside
category with one object $(\Delta)=(sts)=(tst)$ (this name
for the object is natural: the category is isomorphic to its $1$-divided
category).
The atom graph is as follows:
$$
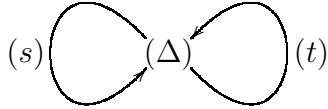
\xy
 (0,0) *++{(\Delta)},
 (-19,0) *++{(s)},
 (19,0) *++{(t)},
(-3,2)="1",
(-3,-2)="2",
(3,-2)="3",
(3,2)="4",
"1";"2" **\crv{(-20,20)&(-20,-20)} ?>*\dir{>}, 
"3";"4" **\crv{(20,-20)&(20,20)} ?>*\dir{>}
\endxy$$
It admits an automorphism of order $2$.

The atom graph of the $3$-divided category is as follows:

$$
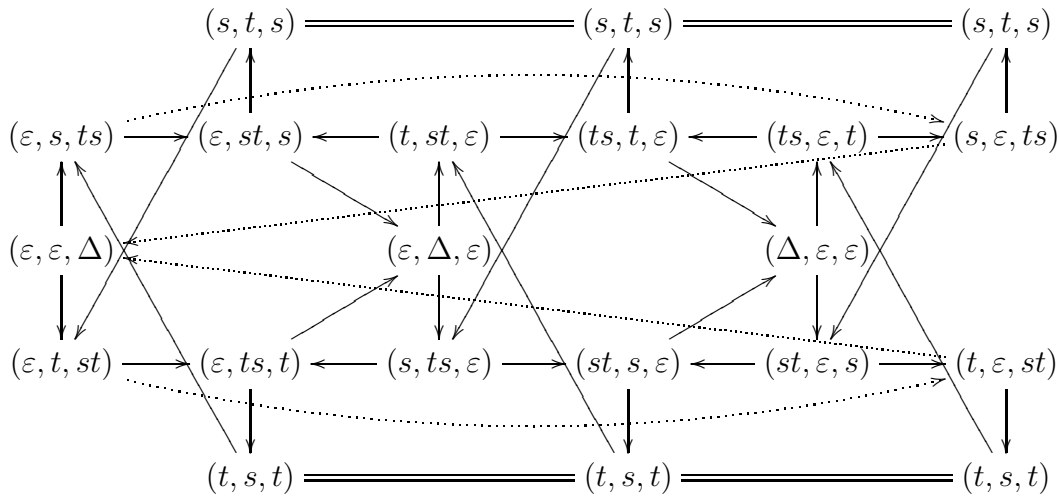
\xymatrix{
& (s,t,s) \ar@{=}[rr] \ar[dddl] & & (s,t,s) \ar@{=}[rr] \ar[dddl]
& & (s,t,s)  \ar[dddl] \\
(\varepsilon,s,ts) \ar[r] \ar@{.>}@/^2em/[rrrrr]
& (\varepsilon,st,s) \ar[u] \ar[dr] &
(t,st,\varepsilon) \ar[r] \ar[l] & (ts,t,\varepsilon) \ar[u] \ar[dr] &
(ts,\varepsilon,t) \ar[r] \ar[l] & (s,\varepsilon,ts) \ar[u]
\ar@{.>}[llllld]\\
(\varepsilon,\varepsilon,\Delta) \ar[u] \ar[d] &  &
(\varepsilon,\Delta,\varepsilon) \ar[u] \ar[d]&  &
(\Delta,\varepsilon,\varepsilon)\ar[u] \ar[d] &  \\
(\varepsilon,t,st) \ar[r]  \ar@{.>}@/_2em/[rrrrr]
& (\varepsilon,ts,t) \ar[d] \ar[ur] &
(s,ts,\varepsilon) \ar[r] \ar[l] & (st,s,\varepsilon) \ar[d] \ar[ur] &
(st,\varepsilon,s) \ar[r] \ar[l] & (t,\varepsilon,st) \ar[d]
\ar@{.>}[lllllu] \\ 
& (t,s,t) \ar@{=}[rr] \ar[uuul] & & (t,s,t)\ar@{=}[rr] \ar[uuul] 
& & (t,s,t) \ar[uuul]  }$$

To improve readability, copies of the vertices $(s,t,s)$ and $(t,s,t)$
have been introduced. Because the graph is better imagined lying on the
surface of a cylinder, we have used dotted arrows to represent the
``hidden'' edges, going behind the cylinder. The graph admits
a symmetry of order $3$ (rotating the cylinder by one third of a turn).
It also admits a symmetry of order $2$ (reflection with horizontal
axis).
The diagram automorphism has order $6$, and is obtained by composing
these two symmetries.

\section{Example: weak Garside groups vs Garside groups}

Let $(\CC,\phi,\Delta)$ be a categorical Garside structure
with Garside groupoid $\CG$.
Let $x\in\CC$ be an object.
It is tempting to think that the category $\CC_x:=\Hom_{\CC}(x,x)$ is 
a Garside category with group of fractions $\CG_x$, the
structure group at $x$.
This is not true, as shown by the following counterexample.

Let $\CC$ be the category defined as the quotient of the free
category on 
$$
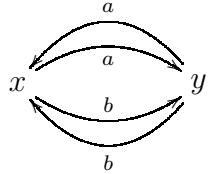
\xymatrix{ x \ar@/^1.2em/[rr]_a \ar@/_1.2em/[rr]^b & & y
\ar@/_2em/[ll]_{a} 
\ar@/^2em/[ll]^{b} }$$
by the relations $a^3=b^3$ (whatever the source may be).
It is a Garside category with Garside
element $\Delta:=a^3=b^3$.
The lattice of left divisors of $\Delta_x^2$ looks as follows:
$$
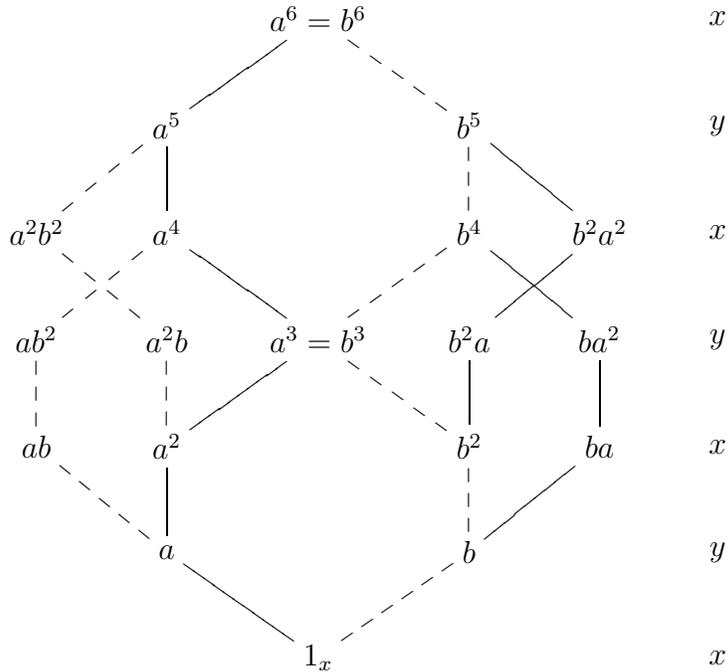
\xymatrix{    & &  a^6=b^6 \ar@{-}[dl] \ar@{--}[dr]  & &  & x\\
& a^5 \ar@{--}[dl] \ar@{-}[d]&  &  b^5  \ar@{--}[d] \ar@{-}[dr]  &  & y\\
a^2b^2\ar@{--}[dr] & a^4  \ar@{-}[dr]  \ar@{--}[dl]&
& b^4 \ar@{--}[dl]  \ar@{-}[dr] & b^2a^2 \ar@{-}[dl] & x \\
ab^2	&	a^2b	& a^3=b^3 	& b^2a 	& ba^2  & y \\
ab\ar@{--}[u] & a^2  \ar@{-}[ur]  \ar@{--}[u]&
& b^2 \ar@{--}[ul]  \ar@{-}[u] & ba \ar@{-}[u]  & x \\
& a \ar@{--}[ul] \ar@{-}[u]&  &  b  \ar@{--}[u] \ar@{-}[ur]  & & y \\
  & & 1_x \ar@{-}[ul] \ar@{--}[ur]  & & & x
}$$
The restriction of the lattice to the $x$ lines is not a lattice:
$a^2$ and $b^2$ do not have a colimit in $\CC_x$, although they
have a colimit in $\CC$ ($a^3=b^3$)
and they have common multiples in $\CC_x$.

\medskip

\section*{Thanks}
Inspiration for this article came as I was working on \cite{kapiun}
and observed a relation between certain topological computations and
a construction of Drew Armstrong, the \emph{$m$-divisible
non-crossing partitions} from \cite{armstrong}.
The connection with Armstrong's work
will be explained in a further version of this preprint.
Although Drew had already told me about his construction,
its importance became clear to me only after
Vic Reiner brought it to my attention, in connection with a joint work
in progress, \cite{br}.
Michel Brou\'e,
Daan Krammer and Jean Michel should also be thanked, for
pleasant discussions at early stages of this work.

\end{document}